\documentclass[a4paper]{article}
\usepackage{amsmath}
\usepackage{amstext}
\usepackage{amsfonts}
\usepackage{euscript}
\usepackage{graphicx}
\usepackage{euscript}
\usepackage[debug]{hyperref}
\usepackage{amssymb}
\input xy
\xyoption{all}

\def\scr{\EuScript}

\newcommand{\N}{\mathbb{N}}
\newcommand{\ZZ}{\mathbb{Z}}
\newcommand{\CC}{\mathbb{C}}
\newcommand{\QQ}{\mathbb{Q}}

\newcommand{\OO}{{\scr O}}

\newcommand{\sbullet}{{\hspace{.1em}\scriptsize\bullet\hspace{.1em}}}

 \DeclareMathOperator{\HS}{HS}

\DeclareMathOperator{\IHS}{IHS} \DeclareMathOperator{\Der}{Der}
\DeclareMathOperator{\Ider}{Ider}

\DeclareMathOperator{\fDer}{\it Der} \DeclareMathOperator{\fIder}{\it Ider}

\DeclareMathOperator{\Aut}{Aut}

\DeclareMathOperator{\diff}{\rm Diff}

\DeclareMathOperator{\supp}{\rm supp}

\newcommand{\pcirc}{{\scriptstyle \,\circ\,}}

\newcommand\D{\underline{D}}

\newcommand\Id{{\rm Id}}

\long\def\inhibe#1\endinhibe{\relax}

\DeclareMathOperator{\Spec}{{\rm Spec}}

\DeclareMathOperator{\Specmax}{{\rm Specmax}}

\newcommand{\dpar}[1]{\partial_{#1}}

\DeclareMathOperator{\ann}{\rm ann}

\renewcommand{\thesubsection}{\arabic{section}.\arabic{subsection}}
\newcounter{numero}[subsection]
\renewcommand{\thenumero}{(\thesubsection .\arabic{numero})}

\newenvironment{corolario}{\medskip
\refstepcounter{numero}\noindent {\sc  \thenumero\ Corollary.}\
\it}{\vspace{1ex}\par}

\newenvironment{teorema}{\medskip
\refstepcounter{numero}\noindent {\sc  \thenumero\ Theorem.}\
\it}{\vspace{1ex}\par}

\newenvironment{lema}{\medskip
\refstepcounter{numero}\noindent {\sc  \thenumero\ Lemma.}\
\it}{\vspace{1ex}\par}

\newenvironment{definicion}{\medskip
\refstepcounter{numero}\noindent {\sc  \thenumero\ Definition.}\
\it}{\vspace{1ex}\par}

\newenvironment{proposicion}{\medskip
\refstepcounter{numero}\noindent {\sc  \thenumero\ Proposition.}\
\it}{\vspace{1ex}\par}

\newenvironment{nota}{\medskip
\refstepcounter{numero}\noindent {\sc  \thenumero\ Remark.}\
}{\vspace{1ex}\par}

\newenvironment{ejemplo}{\medskip
\refstepcounter{numero}\noindent {\sc  \thenumero\ Example.}\
}{\vspace{1ex}\par}

\newenvironment{prueba}{
\noindent {\sc  Proof.}\ }{\hfill Q.E.D.\vspace{1ex}\par}

\newenvironment{question}{\medskip
\refstepcounter{numero}\noindent {\sc  \thenumero\ Question.}\
}{\vspace{1ex}\par}

\newcommand{\numero}{\refstepcounter{numero}\noindent {\sc  \thenumero\ }}

\title{On the modules of $m$-integrable derivations in non-zero characteristic}

\author{Luis Narv\'{a}ez Macarro\thanks{Partially supported by MTM2007-66929, MTM2010-19298
 and FEDER. }}

\date{}

\begin{document}

\maketitle

\begin{abstract} Let $k$ be a commutative ring and $A$ a commutative $k$-algebra.
Given a positive integer $m$, or $m=\infty$,
we say that a $k$-linear derivation $\delta$ of $A$ is $m$-integrable if
it extends up to a Hasse--Schmidt derivation $D=(\Id,D_1=\delta,D_2,\dots,D_m)$ of
$A$ over $k$ of length $m$. This condition is automatically satisfied for any $m$ under one of the following orthogonal hypotheses: (1) $k$ contains the rational
numbers and $A$ is arbitrary, since we can take $D_i=\frac{\delta^i}{i!}$; (2) $k$ is arbitrary and $A$ is a smooth $k$-algebra.
The set of $m$-integrable derivations of $A$ over $k$ is an $A$-module which will be denoted by $\Ider_k(A;m)$. In this paper we prove that, if $A$ is a
finitely presented $k$-algebra and $m$ is a positive integer, then a $k$-linear derivation $\delta$ of $A$ is $m$-integrable if and only if the induced derivation $\delta_{\mathfrak{p}}:A_{\mathfrak{p}} \to A_{\mathfrak{p}} $ is $m$-integrable for each prime ideal $\mathfrak{p}\subset A$. In particular, for any locally
finitely presented morphism of schemes $f:X \to S$ and any positive integer $m$, the $S$-derivations of $X$ which
are locally $m$-integrable form a quasi-coherent submodule $\fIder_S(\OO_X;m)\subset
\fDer_S(\OO_X)$ such that, for any affine open sets $U=\Spec A \subset X$ and
$V=\Spec k \subset S$, with $f(U)\subset V$, we have $\Gamma(U,\fIder_S(\OO_X;m))
=\Ider_k(A;m)$ and $\fIder_S(\OO_X;m)_p = \Ider_{\OO_{S,f(p)}}(\OO_{X,p};m)$ for
each $p\in X$.  We also give, for each positive integer $m$, an algorithm
to decide whether all derivations are $m$-integrable or not.
\medskip

\end{abstract}

\section*{Introduction}

Let us start by recalling the algebraic interpretation of the integration of a vector field.
Let $X$ be a complex algebraic variety and $\chi$ an algebraic vector field on $X$,
or, equivalently, a $\CC$-derivation $\delta:\OO_X\to \OO_X$ of the sheaf of regular
functions. Let us denote by $X[t]= \mathbb{A}^1_{\CC} \times X$, $\CC[\varepsilon] =
\CC[t]/(t^2)$, $X[\varepsilon] = \Spec \CC[\varepsilon] \times X$ and
$\overline{\delta}:X[\varepsilon] \to X$ the map of schemes determined by (and
determining) $\delta$: any section $f$ of $\OO_X$ is mapped to the section $f +
\delta(f)\varepsilon$ of $\OO_X[\varepsilon]$.
\medskip

 If $X$ is
nonsingular, we can consider the flow $\Theta:{\cal U} \to X^{\rm an}$ associated
with $\chi^{\rm an}$, where ${\cal U}\subset X[t]^{\rm an}=\CC \times X^{\rm an}$ is an open neighbourhood
of ${\cal X}=\{0\}\times X^{\rm an}$. It turns out that for any holomorphic (or algebraic)
function $f$ on an open set $V\subset X^{\rm an}$, the function $\Theta^*(f)=f\pcirc \Theta$ is given by
$$ (t,p)\in \Theta^{-1}(V)\subset \CC \times X^{\rm an} \mapsto \sum_{i=0}^\infty t^i \frac{\delta^i(f)}{i!}(p)\in \CC$$
for $|t|$ small enough. Hence, the formal
completion of $\Theta$ along $\cal X$, $\widehat{\Theta}:\widehat{\cal
U}=\widehat{X[t]^{\rm an}} \to X^{\rm an}$, comes from the purely (formal) algebraic map $\widehat{X[t]} \to X$ associated with the 
{\em exponential map} $e^{t\delta}:\OO_X \to \OO_X[[t]]$ attached to $\chi$
(or to $\delta$) defined as
$$e^{t\delta}(f) = \sum_{i=0}^\infty t^i \frac{\delta^i(f)}{i!}$$
for any regular function $f$ on some Zariski open set of $X$.
\medskip

The exponential map
$e^{t\delta}$ is a lifting of $\delta$ (it coincides with $\overline{\delta}$ $\mod
t^2$) and it can be regarded as the algebraic incarnation of the integration of the
vector field $\chi$.\medskip

The exponential map of a vector field makes sense not only over the complex numbers,
but over any field of characteristic zero, and in fact it also works if $X$ is
eventually singular. However, it does not make sense over a field $k$ of positive
characteristic. \medskip

Nevertheless, the notion of {\em Hasse--Schmidt derivation} allows
us to define what integrability means for a vector field in such a case (see
\cite{brown_1978,mat-intder-I}). Given a commutative ring $k$ and a commutative $k$-algebra $A$, a Hasse--Schmidt derivation of $A$ over $k$ (of length $\infty$) is a sequence
$D=(\Id,D_1,D_2,D_3,\dots)$ of $k$-linear operators of $A$ which appear as the coefficients of a $k$-algebra map $\Phi:A \to A[[t]]$ such that $\Phi(a) \equiv a \mod t$ for all $a\in A$:  $\Phi(a) = a + D_1(a) t + D_2(a) t^2 + \cdots$. That property is equivalent to the fact that the $D_i$ satisfy the Leibniz equality:
$$ D_0 =\Id,\quad D_i(ab) = \sum_{r+s=i} D_r(a) D_s(b)\quad \forall a,b\in A,\ \forall i\geq 1.$$
A $k$-linear derivation $\delta:A\to A$ is said to be ($\infty$-)integrable if there is a Hasse--Schmidt derivation $D$ of $A$ over $k$ (of length $\infty$) such that $D_1=\delta$, or in other words, if the $k$-algebra map $\overline{\delta}: a\in A \mapsto a + \delta(a) \varepsilon \in A[\varepsilon]= A[[t]]/(t^2)$ can be lifted up to a $k$-algebra map $\Phi: A \to A[[t]]$.
The set of $k$-linear derivations of $A$ which are integrable is a submodule of $\Der_k(A)$, which is denoted by $\Ider_k(A)$.
\medskip

When $A$ is a smooth $k$-algebra over an arbitrary commutative ring $k$ or when $k$ contains the rational numbers, any $k$-linear derivation $\delta:A\to A$ is ($\infty$-)integrable.
The modules $\Ider_k(A)$, and more
generally, the Hasse--Schmidt derivations of $A$ over $k$ seem to play an important role
among the differential structures in Commutative Algebra and Algebraic Geometry
(see \cite{vojta_HS}, \cite{nar_2009}). They behave better in positive characteristic
than $\Der_k(A)$ (see for instance \cite{molinelli_1979} or \cite{seiden_1966}) and one expects that they can help to understand (some of) the differences between
singularities in zero and non-zero characteristics, but they are difficult to deal with. For instance, it is not clear at all that ($\infty$-)integrability is a local property (in the sense that can be tested locally at the primes ideals of $A$).
\medskip

For a given positive integer $m$, the $m$-integrability of a $k$-linear derivation $\delta:A\to A$ is defined as the existence of a $k$-algebra map $\Phi: A \to A[[t]]/(t^{m+1})$ lifting the map $\overline{\delta}$ defined above. The set of $k$-linear derivations of $A$ which are $m$-integrable is a submodule of $\Der_k(A)$, which is denoted by $\Ider_k(A;m)$. One obviously has
$\Der_k(A)=\Ider_k(A;1) \supset \Ider_k(A;2) \supset \Ider_k(A;3) \supset \cdots \supset \Ider_k(A;\infty)=\Ider_k(A)$.\medskip

This paper is devoted to the study of the modules $\Ider_k(A;m)$, for $m\geq 1$.
\medskip

One of the main difficulties when dealing with
$m$-integrability of a derivation is that one cannot proceed step by step: a
derivation $\delta$ can be $(m+r)$-integrable, but it may have an intermediate $m$-integral $D=(\Id,D_1=\delta,D_2,\dots,D_m)$
which does not  extends up to a Hasse--Schmidt derivation of length $(n+r)$ (cf. Example 3.7 in \cite{nar_2009}).
\smallskip

Our main results are the following:
\smallskip

\noindent (I) If $A$ is a finitely presented $k$-algebra and $m$ is a positive integer, then the property of being
$m$-integrable for a $k$-derivation $\delta$ of $A$ is a local property, i.e. $\delta$ is $m$-integrable if and only if the induced derivation $\delta_{\mathfrak{p}}:A_{\mathfrak{p}} \to A_{\mathfrak{p}} $ is $m$-integrable for each prime ideal $\mathfrak{p}\subset A$. As a consequence, for any locally
finitely presented morphism of schemes $f:X \to S$ and any positive integer $m$, the $S$-derivations of $X$ which
are locally $m$-integrable form a quasi-coherent submodule $\fIder_S(\OO_X;m)\subset
\fDer_S(\OO_X)$ such that, for any affine open sets $U=\Spec A \subset X$ and
$V=\Spec k \subset S$, with $f(U)\subset V$, we have $\Gamma(U,\fIder_S(\OO_X;m))
=\Ider_k(A;m)$ and $\fIder_S(\OO_X;m)_p = \Ider_{\OO_{S,f(p)}}(\OO_{X,p};m)$ for
each $p\in X$ (see Theorem \ref{teo:criter-loc-ider-fp}
and Corollary \ref{cor:Ider_fp_schemes}). We have then a decreasing sequence of quasi-coherent modules
$$\fDer_S(\OO_X)= \fIder_S(\OO_X;1) \supset \fIder_S(\OO_X;2) \supset \fIder_S(\OO_X;3)\supset \cdots$$
and all the quotients $\fDer_S(\OO_X)/\fIder_S(\OO_X;m)$ are supported by the non-smooth\-ness locus of $f:X\to S$.
\smallskip

\noindent (II) For a given $k$-algebra $A$ and for any positive integer $m$, there
is a constructive procedure to see whether \underline{all} $k$-derivations of $A$
are $m$-integrable or not. In particular, if $A$ and $k$ are ``computable'' rings,
then the above procedure becomes an effective algorithm (although of exponential
complexity with respect to $m$) to decide whether the equality
$\Ider_k(A;m)=\Der_k(A)$ is true or not (see \ref{subsect:algo}). \smallskip

Let us now comment on the content of this paper. \smallskip

In section 1 we review the notion of Hasse--Schmidt derivation and its basic
properties. We study logarithmic Hasse--Schmidt derivations with respect to an ideal
$I$ of some ambient algebra $A$ and their relationship with Hasse--Schmidt
derivations of the quotient $A/I$. In the last part we focus on the description of
Hasse--Schmidt derivations on polynomial or power series algebras.
\smallskip

Section 2 contains the main results of this paper. First, we define
$m$-integrability and logarithmic $m$-integrability and give a characterization of
$(m+1)$-integrability for a Hasse--Schmidt derivation of length $m$. In section
\ref{subsect:jacob} we give some criteria for a derivation to be integrable, based
on and extending previous results of \cite{mat-intder-I} and \cite{traves-2000}.
Next, we study the behaviour of $m$-integrability under localization, for finite
$m$, and prove (I) above. In the last part we prove the results needed to justify
procedure (II) above. \smallskip

In Section 3 we first compute some concrete examples and illustrate the nonlinear
equations one encounters when computing systems of generators of the modules
$\Ider_k(A;m)$. In the second part we state some questions, which seem to be
important for understanding the relationship between the modules of $m$-integrable
derivations and singularities.
\medskip

I would like to thank Herwig Hauser and Orlando Villamayor for many helpful and
inspiring discussions, and also Herwig Hauser for proposing the last two examples in
section 3. I would also like to thank Herwig Hauser and Eleonore Faber for some comments and suggestions on a previous version of this paper.

\section{Notations and preliminaries}

\subsection{Notations}

Throughout the paper we will use the following notations:
\smallskip

\noindent -) $k$ will be a commutative ring and $A$ a commutative
$k$-algebra.\smallskip

\noindent -) $\N_+ := \{n\in\N\ |\ n\geq 1\}$, $\overline{\N} := \N \cup
\{\infty\}$, $\overline{\N}_+ := \N_+ \cup \{\infty\}$.\smallskip

\noindent -) If $n\in \N_+$, $[n]:=\{0,1,\dots,n\}$, $[n]_+:=[n]\cap \N_+$ and $[\infty]:=\N$.
\smallskip

\noindent -) If $n\in \N_+$, $A_n :=A[[t]]/(t^{n+1})$ and $A_\infty
=A[[t]]$.  Each $A_n$ is an augmented $A$-algebra, the
augmentation ideal  $\ker (A_n \to A)$ being generated by $t$.\smallskip

\noindent -) For $n\in \overline{\N}_+$ and $m\in [n]_+$, let us denote
by $\pi_{nm}: A_n \to A_m$ the natural epimorphism of augmented
$A$-algebras.\smallskip

\noindent -) If $\alpha =(\alpha_1,\dots,\alpha_d)\in\N^d$,  $\supp \alpha = \{r\in\{1,\dots,d\}\ |\ \alpha_r\neq 0\}$ and $|\alpha|:= \alpha_1+\cdots+\alpha_d$.
\smallskip

\noindent -) The ring of $k$-linear differential operators of $A$ will be denoted by
$\diff_{A/k}$ (see \cite{ega_iv_4}).\smallskip

\noindent -) For $A=k[x_1,\dots,x_d]$ or $A=k[[x_1,\dots,x_d]]$, we will denote by
$\partial_r:A\to A$ the partial derivative with respect to $x_r$.

\subsection{Hasse-Schmidt derivations}

In this section we remind the definition and basic facts of Hasse--Schmidt
derivations (see \cite{has37},\cite{mat_86}, \S 27, and \cite{traves-phd},
\cite{vojta_HS}, \cite{nar_2009} for more recent references). We also introduce the
basic constructions that will be used throughout the paper.

\begin{definicion} A {\em Hasse--Schmidt derivation} of $A$ (over $k$)  of length $n\geq 1$
(resp. of length $\infty$) is a sequence $D=(D_i)_{i\in [n]}$
 of $k$-linear maps $D_i:A
\longrightarrow A$, satisfying the conditions: $$ D_0=\Id_A, \quad
D_i(xy)=\sum_{r+s=i}D_r(x)D_s(y) $$ for all $x,y \in A$ and for all
$i\in [n]$. We denote by
$\HS_k(A;n)$ the set of all Hasse--Schmidt derivations of $A$ (over
$k$) of length $n\in \overline{\N}$ and $\HS_k(A)=\HS_k(A;\infty)$.
\end{definicion}

\noindent \numero The $D_1$ component of any Hasse-Schmidt derivation $D\in\HS_k(A;n)$
is a $k$-derivation of $A$. More generally, the $D_i$ component is a $k$-linear differential operator of order $\leq i$ with $D_i(1)=0$ for $i=1,\dots,n$.
\smallskip

\noindent\numero Any Hasse--Schmidt derivation $D\in\HS_k(A;n)$ is determined by the
$k$-algebra homomorphism $\Phi_D: A \to A_n$ defined by $\Phi_D(a) = \sum_{i=0}^n
D_i(a)t^i$ and satisfying $\Phi_D(a)\equiv a \mod t$. The $k$-algebra homomorphism
$\Phi_D$ can be uniquely extended to a $k$-algebra automorphism $\widetilde{\Phi}_D:
A_n \to A_n$ with $\widetilde{\Phi}_D(t)=t$:
$$ \widetilde{\Phi}_D\left(\sum_{i=0}^n a_i t^i\right) = \sum_{i=0}^n \Phi(a_i) t^i.$$
So, there is a bijection between $\HS_k(A;n)$ and the subgroup of
$\Aut_{k-\text{alg}}(A_n)$ consisting of the automorphisms
$\widetilde{\Phi}$ satisfying $\widetilde{\Phi}(a) \equiv a \mod t$
for all $a\in A$ and $\widetilde{\Phi}(t)=t$. In particular,
$\HS_k(A;n)$ inherits a canonical group structure which is
explicitly given by $D\pcirc D' = D''$ with $ D''_{l} = \sum_{i+j=l}
D_i \pcirc D'_j$, the identity element of $\HS_k(A;n)$ being
$(\Id_A,0,0,\dots)$. It is clear that the map
$(Id_A,D_1)\in \HS_k(A;1) \mapsto D_1 \in \Der_k(A)$
is
an isomorphism of groups, where we consider the addition as internal operation in $\Der_k(A)$.
\smallskip

\noindent \numero For any $a\in A$ and any $D\in\HS_k(A;n)$, the sequence $a\sbullet
D$ defined by $(a\sbullet D)_i=a^i D_i$, $i\in [n]$, is again a Hasse--Schmidt
derivation of $A$ over $k$ of length $n$ and $\Phi_{a\sbullet D}(b)(t) =
\Phi_D(b)(at)$ for all $b\in A$. We have $(a a')\sbullet D= a\sbullet (a'\sbullet
D)$, $1\sbullet D=D$ and $0\sbullet D=$ the identity element.
\smallskip

\noindent\numero \label{nume:truncation-inverse-limit} For $1\leq m \leq n\in
\overline{\N}$, let us denote by $\tau_{nm}: \HS_k(A;n) \to \HS_k(A;m)$ the {\em
truncation} map defined in the obvious way. One has $\Phi_{\tau_{nm} D} =
\pi_{nm}\pcirc \Phi_D$. Truncation maps are group homomorphisms and they satisfy
$\tau_{nm}(a\sbullet D)=a\sbullet \tau_{nm}D$. It is clear that the group
$\displaystyle \HS_k(A;\infty)$ is the inverse limit of the groups $\HS_k(A;m)$,
$m\in\N$.

\begin{definicion} Let $q\geq 1$ be an integer or $q=\infty$, and $D\in\HS_k(A;q)$.
For each integer $m\geq 1$ we define $D[m]$ as the Hasse--Schmidt derivation (over
$k$) of length $mq$ determined by the $k$-algebra map obtained by composing the
following maps:
$$ A \xrightarrow{\Phi_D} A_q = A[[t]]/(t^{q+1})
\xrightarrow{\overline{t} \mapsto \overline{t}^m} A_{mq}= A[[t]]/(t^{mq+1}).$$
In the case $q=1$ and $D=(\Id_A,\delta)$, we simply denote $\delta[m]:=D[m]$.
\end{definicion}

If $D=(\Id_A,D_1,D_2,\dots)\in\HS_k(A;q)$, then
$$D[m] =(\Id_A,0,\dots,0,\stackrel{\scriptsize\underbrace{m}}{D_1},0,\dots,0,
\stackrel{\scriptsize\underbrace{2m}}{D_2},0,\dots)\in \HS_k(A;mq).$$ The map $ D\in
\HS_k(A;q) \mapsto D[m] \in \HS_k(A;qm)$ is a group homomorphism and we have
$(a^m\bullet D)[m] = a\bullet D[m]$, $(\tau_{qq'}D)[m] = \tau_{qm,q'm}(D[m])$ for
$a\in A, 1\leq q'\leq q$.

\begin{definicion} \label{def:ell} For each $n\in\overline{\N}_+$ and each
$E\in\HS_k(A;n)$, we denote $\ell(E)=0$ if $E_1\neq 0$, $\ell(E)=n$ if $E$ is the
identity and $\ell(E)=$ maximun of the $r\in [n]$ such that $E_1=\cdots=E_r=0$
otherwise.
\end{definicion}

\begin{definicion} Let $I\subset A$ be an ideal and $m\in\overline{\N}_+$. We say that:
\begin{enumerate}
\item[1)] A $k$-derivation $\delta:A\to A$ is {\em $I$-logarithmic} if $\delta(I)\subset I$. The set of
$k$-linear derivations of $A$ which are $I$-logarithmic is
denoted by $\Der_k(\log I)$.
\item[2)]  A Hasse--Schmidt derivation $D\in\HS_k(A;m)$ is called {\em $I$-logarithmic} if $D_i(I)\subset
I$ for any $i\in [m]$. The set of Hasse--Schmidt derivations $D\in\HS_k(A;m)$ which are $I$-logarithmic is denoted by $\HS_k(\log I;m)$. When $m=\infty$ it will be simply denoted by $\HS_k(\log I)$.
\end{enumerate}
\end{definicion}

The set $\Der_k(\log I)$ is obviously a $A$-submodule of $\Der_k(A)$. Any $\delta\in
\Der_k(\log I)$ gives rise to a unique $\overline{\delta}\in \Der_k(A/I)$ satisfying
$\overline{\delta}\pcirc \pi = \pi \pcirc \delta$, where $\pi:A\to A/I$ is the
natural projection. Moreover, if
$A=k[x_1,\dots,x_d]$ or $A=k[[x_1,\dots,x_d]]$, the sequence of $A$-modules
\begin{equation*}
0 \to I\Der_k(A) \xrightarrow{\text{incl.}} \Der_k(\log I)
\xrightarrow{\delta\mapsto \overline{\delta}} \Der_k(A/I)\to 0
\end{equation*}
is exact. \smallskip

\numero \label{num:HS-induced-by-HS-log} In the same vein, the set $\HS_k(\log I;m)$
is a subgroup of $\HS_k(A;m)$ and we have $A\sbullet \HS_k(\log I;m)\subset
\HS_k(\log I;m)$, $\HS_k(\log I;m)[n]\subset \HS_k(\log I;mn)$, $n\in\N$. A $D\in
\HS_k(A;m)$ is $I$-logarithmic if and only if its corresponding $k$-algebra
homomorphism $\Phi_D:A\to A_m$ satisfies $\Phi_D(I) \subset I_m:= \ker \pi_m$, where
$\pi_m:A_m \to \left(A/I\right)_m$ is the natural projection\footnote{Observe that
$\ker \pi_m = I A_m$ when $I$ is finitely generated or $m$ is finite.}. Moreover, a
$I$-logarithmic Hasse--Schmidt derivation $D\in \HS_k(\log I;m)$ gives rise to a
unique $\overline{D}\in \HS_k(A/I;m)$ such that $\overline{D}_i\pcirc \pi = \pi
\pcirc D_i$ for all $i\in [m]$, and the following diagram is commutative
\begin{equation*}
\xymatrix{ A \ar[r]^{\Phi_D} \ar[d]^{\pi} & A_m \ar[d]^{\pi_m} \\
A/I  \ar[r]^{\Phi_{\overline{D}}} & \left(A/I\right)_m.}
\end{equation*}
The map $\Pi_m: D\in \HS_k(\log I;m) \to \overline{D} \in \HS_k(A/I;m)$ is clearly a
homomorphism of groups and $\Pi_m(a\sbullet D)= \pi(a)\sbullet \Pi_m(D)$. So, its
kernel contains the subgroup $I\sbullet \HS_k(A;m)$ generated by the $a\sbullet E$,
with $a\in I$ and $E\in\HS_k(A;m)$. It is also clear that $\tau_{mn}\pcirc \Pi_m =
\Pi_n \pcirc \tau_{mn}$ and $(\Pi_m D)[n]= \Pi_{mn}(D[n])$. \smallskip

\noindent \numero  \label{nume:basic-localiz}
 Let $S\subset A$ be a multiplicative set. For each $k$-linear
differential operator $P:A\to A$, let us denote by $\widetilde{P}:S^{-1}A \to
S^{-1}A$ its canonical extension. We know that the map $P\in\diff_{A/k} \mapsto
\widetilde{P}\in\diff_{S^{-1}A/k}$ is a ring homomorphism. Let $m\geq 1$ be an
integer or $m=\infty$ and $\mathfrak{a}\subset A$ an ideal. Here is a summary of the
basic facts of the behaviour of Hasse-Schmidt derivations under localization:
\smallskip

\noindent -) For any $D=(D_i)\in\HS_k(A;m)$, the sequence $\widetilde{D}:=
(\widetilde{D_i})$ is a Hasse-Schmidt derivation of $S^{-1}A$ (over $k$ of length
$m$) and the following diagram is commutative
\begin{equation*}
\xymatrix{ A \ar[r]^{\Phi_D} \ar[d]^{\text{can.}} & A_m \ar[d]^{\text{can.}} \\
S^{-1}A  \ar[r]^{\Phi_{\widetilde{D}}} & (S^{-1}A)_m.}
\end{equation*}
Moreover, if $D$ is $\mathfrak{a}$-logarithmic, then $\widetilde{D}$ is
$(S^{-1}\mathfrak{a})$-logarithmic.\smallskip

\noindent -) The map $\Theta_m:D\in \HS_k(A;m) \to \widetilde{D}\in \HS_k(S^{-1}
A;m)$ is a group homomorphism, $\Theta_m(a\sbullet D) = \frac{a}{1}
\sbullet\Theta_m(D)$ and the following diagram is commutative:
\begin{equation*}
\xymatrix{ \HS_k(\log \mathfrak{a};m) \ar[r]^{\Theta_m} \ar[d]^{\Pi_m} & \HS_k(\log (S^{-1}\mathfrak{a});m) \ar[d]^{\Pi_m} \\
\HS_k(A/\mathfrak{a};m)  \ar[r]^{\Theta_m} & \HS_k(S^{-1}A/S^{-1}\mathfrak{a};m).}
\end{equation*}
Moreover, $\tau_{mn}\pcirc \Theta_m = \Theta_n \pcirc \tau_{mn}$ and $(\Theta_m
D)[n]= \Theta_{mn}(D[n])$. \smallskip

The extension of Hasse--Schmidt derivations to rings of fractions is a particular
case of the formally \'etale extensions (cf. \cite{masson_PhD} and \cite{traves-2000}, th. 1.5).

\subsection{Hasse--Schmidt derivations of polynomial or formal power series algebras}

Throughout this section we assume that $A=k[x_1,\dots,x_d]$ or $A=k[[x_1,\dots,x_d]]$.
The {\em Taylor differential operators} $\Delta^{(\alpha)}:A\to A$, $\alpha\in\N^d$,
are defined by:
$$ g(x_1+T_1,\dots,x_d+T_d) =\sum \Delta^{(\alpha)}(g)
T^\alpha,\quad \forall g\in A.$$ It is well known that
$\{\Delta^{(\alpha)}\}_{|\alpha|\leq i}$ is a basis of the left (resp. right)
$A$-module of $k$-linear differential operators of $A$ of order $\leq i$. So, if
$D\in \HS_k(A;m)$, there are unique $C^i_{\alpha}\in A$, $\alpha\in\N^d$, $0<i\leq
|\alpha|\in [m]_+$, such that $D_i = \sum_{0<|\alpha|\leq i} C^i_{\alpha}
\Delta^{(\alpha)}$, $i\in [m]_+$. On the other hand, there are unique $c_{ri}\in A$,
$i\in [m]_+$, $1\leq r\leq d$, such that
$$\Phi_D(x_r)= x_r + \sum_{i=1}^m   c_{ri} t^i,\quad 1\leq r\leq d.$$
In fact, any system of $c_{ri}\in A$, $i\in [m]_+$, $1\leq r\leq d$, determines uniquely such a homomorphism of $k$-algebras $A\to A_m$ and so a Hasse--Schmidt derivation $D\in \HS_k(A;m)$.
\medskip

The following proposition gives the relationship between the $C^i_{\alpha}$ and the
$c_{ri}$ above. Its proof does not contain any surprise and it is left up to the
reader.

\begin{proposicion} \label{prop:formulon} With the above notations, the following properties hold:
\begin{enumerate}
\item[1)] $c_{ri}= D_i(x_r) = C^i_{e_r}$, with $e_r = (0,\dots,\stackrel{\scriptsize\underbrace{r}}{1},\dots,0)$, for all $i\in [m]_+$, $r=1,\dots,d$.
\item[2)] $$C_\alpha^i=\sum_{\substack{\{\varepsilon_r\}_{r\in\supp \alpha}\\ \varepsilon_r\geq \alpha_r,
|\varepsilon|=i}} \left(
 \prod_{r\in \supp \alpha} \left(
 \sum_{\substack{\beta_1+\cdots+\beta_{\alpha_r}=\varepsilon_r\\ \beta_k> 0}} \prod_{k=1}^{\alpha_r} c_{r,\beta_k}\right)\right)$$
 for all $\alpha\in\N^d$, $|\alpha|\in [m]_+$, $0<i\leq |\alpha|$.
\end{enumerate}
\end{proposicion}

The above proposition is a particular case of Theorem 2.8 in \cite{magda_nar_hs}.
For the sake of completeness we include, without proof, the following result.

\begin{proposicion} Let $C^i_{\alpha}\in A$, $\alpha\in\N^d$, $0<i\leq |\alpha|\in [m]_+$, be a system of elements of $A$
and define $D_0=\Id_A$, $D_i = \sum_{0<|\alpha|\leq i} C^i_{\alpha} \Delta^{(\alpha)}$, $i\in [m]_+$. The following properties are equivalent:
\begin{enumerate}
\item[(a)] The sequence $D=(D_i)_{i\in [m]}$ is a Hasse--Schmidt derivation of $A$ over $k$ of length $m$.
\item[(b)] For all $i\in [m], i\geq 2$, for all $\varrho\in\N^d$ with $2\leq |\varrho|\leq i$ and for all $\beta,\gamma\in\N^d$ with
$\varrho=\beta+\gamma$, $|\beta|,|\gamma| >0$ we have
$\displaystyle C^i_{\varrho}= \sum C^j_{\beta} C^l_{\gamma}$, where the summation indexes are the $(j,l)$ with $j\geq |\beta|, l\geq |\gamma|$ and $j+l=i$.
\end{enumerate}
\end{proposicion}

Let us notice that, if the equivalent properties of the preceding proposition hold, then the $C^i_{\alpha}$ with $2\leq
|\alpha|\leq i$ are determined by the $C^j_{\beta}$ with $1\leq |\beta|\leq j\leq i-1$. This applies in particular to the symbol of the $D_i$,
$\sigma(D_i)=\sum_{|\alpha|=i} C^i_{\alpha} \xi^{\alpha}$, which only depend on $D_1$ (compare with Proposition 2.6 in \cite{nar_2009}).

\begin{definicion} \label{def:Taylor-HS}
The {\em Taylor Hasse-Schmidt derivations} of $A$
 are the $$\Delta^{(s)} := (\Id_A,\Delta_1^{(s)},\Delta_2^{(s)},
\Delta_3^{(s)},\dots ) \in \HS_k(A),\quad  1\leq s\leq d,$$ where $\Delta^{(s)}_i =
 \Delta^{(0,\dots,\stackrel{\scriptsize\underbrace{s}}{i},\dots,0)}$ for each $i\geq 1$.
\end{definicion}

\begin{proposicion}  \label{prop:surjec-HS-log} Assume that $R=k[x_1,\dots,x_d]$,
$S\subset R$ is a multiplicative set and $A=S^{-1}R$ or $A=k[[x_1,\dots,x_d]]$. For
any ideal $I\subset A$, the group homomorphisms $\Pi_m: \HS_k(\log I;m) \to
\HS_k(A/I;m)$, $m\in\overline{\N}$, (see \ref{num:HS-induced-by-HS-log}) are
surjective.
\end{proposicion}

\begin{prueba} Let us prove the proposition in the case $A=S^{-1}R$, the case
$A=k[[x_1,\dots,x_d]]$ being completely similar. Let us call $\sigma:R\to A, \pi:A
\to A/I, \pi_m: A_m \to (A/I)_m$ the canonical maps and let $E\in \HS_k(A/I;m)$ be
any Hasse--Schmidt derivation. Let $a_{ri}\in A$ be elements such that
$$ \Phi_E(\pi(\sigma(x_r))) = \pi(\sigma(x_r)) +\sum_{i\in [m]} \pi(a_{ri}) t^i \in (A/I)_m,
\quad r=1,\dots,d,$$
and let $\Psi:R \to A_m$ be the $k$-algebra map defined by
$$ \Psi(x_r) = \sigma(x_r) +\sum_{i\in [m]} a_{ri} t^i \in A_m,\quad r=1,\dots,d.$$
Since $\Psi(f)\equiv \sigma(f) \mod t$ for each $f\in R$, we deduce that $\Psi(s)$
is invertible for all $s\in S$ and the map $\Psi$ induces $\widetilde{\Psi}:A\to
A_m$. It is clear that $\widetilde{\Psi}(a)\equiv a \mod t$ for each $a\in A$ and
$\pi_m\pcirc \widetilde{\Psi} = \Phi_E\pcirc \pi$. So, $\widetilde{\Psi}$ induces a
$I$-logarithmic Hasse-Schmidt derivation $D\in\HS_k(\log I;m)$ such that
$\Pi_m(D)=E$ (see \ref{num:HS-induced-by-HS-log}).
\end{prueba}

\begin{proposicion}  \label{prop:surjec-HS-log-localiz-poly}
Assume that $R=k[x_1,\dots,x_d]$, $S\subset R$ is a multiplicative set and let
$\mathfrak{a}\subset R$ be a finitely generated ideal. For any (finite) integer
$m\geq 1$, the map
$$ (s,D)\in S \times \HS_k(\log \mathfrak{a};m)\mapsto \frac{1}{s}\sbullet
\Theta_m(D) \in\HS_k(\log (S^{-1}\mathfrak{a});m)$$ is surjective.
\end{proposicion}

\begin{prueba} Let
$E\in\HS_k(\log (S^{-1}\mathfrak{a});m)$ be any $(S^{-1}\mathfrak{a})$-logarithmic
Hasse--Schmidt derivation. Since $m$ is finite, there are $a_{ij}\in R$, $1\leq i=1
\leq d$, $1\leq j\leq m$ and $\sigma\in S$ such that
$$ \Phi_E\left(\frac{x_i}{1}\right) = \frac{x_i}{1} + \left(\frac{a_{i1}}{\sigma}\right) t +
\cdots + \left(\frac{a_{im}}{\sigma}\right) t^m \in (S^{-1}R)_m,\quad i=1,\dots,d.$$
Let us consider the $k$-algebra map $\Phi^0: R \to R_m$ given by
$$ \Phi^0(x_i) = x_i + a_{i1} t + \sigma a_{i2} t^2 +\cdots + \sigma^{m-1} a_{im} t^m \in
R_m,\quad i=1,\dots,d$$ and the corresponding Hasse-Schmidt derivation $D^0\in
\HS_k(R;m)$ with $\Phi^0 =\Phi_{D^0}$. It is clear that
$\left(\frac{\sigma}{1}\right)\sbullet E = \Theta_m(D^0)$. Let $f_1,\dots,f_u\in
\mathfrak{a}$ be a finite system of generators. Since $\Theta_m(D^0)$ is
$(S^{-1}\mathfrak{a})$-logarithmic, we deduce the existence of a $\tau\in S$ such
that $\tau\Phi_{D^0}(f_l) \in A_m \mathfrak{a}$ for all $l=1,\dots,u$. So,
$D:=\tau\sbullet D^0$ is $\mathfrak{a}$-logarithmic and $E =
\left(\frac{1}{\sigma\tau}\right)\sbullet \Theta_m(D)$.
\end{prueba}

Proposition \ref{prop:surjec-HS-log-localiz-poly} is false for $m=\infty$, as shown
for instance in example 1.4 in \cite{traves-2000}.

\begin{corolario} \label{cor:surjec-HS-log-localiz-fp}
Assume that $A$ is a finitely presented $k$-algebra and let $T\subset A$ be a
multiplicative set. Then, for any (finite) integer $m\geq 1$, the map
$$ (t,E)\in T \times \HS_k(A;m)\mapsto \frac{1}{t}\sbullet
\Theta_m(E) \in\HS_k(T^{-1}A;m)$$ is surjective.
\end{corolario}

\begin{prueba} We may assume that $A=R/\mathfrak{a}$
with $R=k[x_1,\dots,x_d]$ and $\mathfrak{a}\subset R$ a finitely generated ideal.
Denote by $\pi:R\to A$ the natural projection and $S=\pi^{-1}(T)$. We have $T^{-1}A
= S^{-1}R/S^{-1}\mathfrak{a}$. Let us look at the following commutative diagram
\begin{equation*}
\xymatrix{S\times \HS_k(\log \mathfrak{a};m)  \ar[r]^{} \ar[d]^{\pi\times\Pi_m} &
\HS_k(\log (S^{-1}\mathfrak{a});m) \ar[d]^{\Pi_m} \\
T \times \HS_k(A;m)  \ar[r]^{} & \HS_k(T^{-1}A;m).}
\end{equation*}
The vertical arrows are surjective by Proposition \ref{prop:surjec-HS-log}. To
conclude, we apply Proposition \ref{prop:surjec-HS-log-localiz-poly}.
\end{prueba}

\section{Integrability}

\subsection{Integrable Hasse--Schmidt derivations}

In this subsection, $A$ will be again an arbitrary $k$-algebra.

\begin{definicion} (Cf. \cite{brown_1978,mat-intder-I})  \label{def:HS-integ}
 We say that a $k$-derivation $\delta:A\to
A$ is {\em $n$-integrable} (over $k$), with $n\in \overline{\N}$, if there is a
Hasse--Schmidt derivation $D\in \HS_k(A;n)$ such that $D_1=\delta$. A such $D$ will
be called a  {\em $n$-integral} of $\delta$. The set of $n$-integrable
$k$-derivations of $A$ is denoted by $\Ider_k(A;n)$. We simply say that $\delta$ is
{\em integrable} if it is $\infty$-integrable and we denote
$\Ider_k(A)=\Ider_k(A;\infty)$.\\

\noindent More generally, we say that a Hasse--Schmidt derivation $D'\in\HS_k(A;m)$
is {\em $n$-integrable} (over $k$), with $m,n\in \overline{\N}, n\geq m$,  if there
is a Hasse--Schmidt derivation $D\in \HS_k(A;n)$ such that $\tau_{nm}D=D'$. A such
$D$ will be called a {\em $n$-integral} of $D'$. The set of $n$-integrable
Hasse--Schmidt derivations of $A$ over $k$ of length $m$ is denoted by
$\IHS_k(A;m;n)$. We simply say that $D'$ is {\em integrable} if it is
$\infty$-integrable and we denote $\IHS_k(A;m) = \IHS_k(A;m;\infty)$.
\end{definicion}

It is clear that the $\Ider_k(A;n)$ are $A$-submodules of $\Der_k(A)$, $
\Der_k(A)=\Ider_k(A;1)\supset\Ider_k(A;2)\supset \Ider_k(A;3)\supset \cdots$ and
\begin{equation} \label{eq:intersection-ider}
 \Ider_k(A) \subset \bigcap_{n\in\N_+} \Ider_k(A;n).
\end{equation}
  It is also clear that the $\IHS_k(A;m;n)$ are subgroups of $\IHS_k(A;m)$, stable by the $\sbullet$
  operation, $\IHS_k(A;m) = \IHS_k(A;m;m) \supset \IHS_k(A;m;m+1)\supset \cdots$ and
\begin{equation} \label{eq:intersection-iHS}
 \IHS_k(A;m) \subset \bigcap_{n\geq m} \IHS_k(A;m;n).
\end{equation}

\begin{ejemplo} \label{ex:1} (1) Let $n\geq 1$ be an integer. If $n!$ is invertible in $A$,
then any $k$-derivation $\delta$ of $A$ is $n$-integrable: we can
take $D \in \HS_k(A;n)$ defined by $D_i=\frac{\delta^i}{i!}$ for
$i=0,\dots,n$. In the case $n=\infty$, if $\QQ\subset A$, one proves
in a similar way that any $k$-derivation of $A$ is integrable.
\smallskip

\noindent (2) If  $A$ is $0$-smooth (i.e. formally smooth for the discrete topologies)
$k$-algebra, then any $k$-derivation of $A$ is integrable (cf. \cite{mat_86}, Theorem
27.1).
\end{ejemplo}

\begin{nota} \label{nota:taylor-HS} A particularly important case of example \ref{ex:1} is
 $A=k[x_1,\dots,x_d]$ or $A=k[[x_1,\dots,x_d]]$.
 In this case we can do better than in example \ref{ex:1} and even
 exhibit a special integral for each $D\in\HS_k(A;m)$, $m\in\N_+$. Namely, consider the Hasse--Schmidt derivation
 $\varepsilon(D)\in\HS_k(A)$ determined by the $k$-algebra map $A=k[x_1,\dots,x_d] \to A[[t]]$ sending each $x_r$ to
 $ \sum_{i\in [m]} D_i(x_r) t^i \in A[[t]]$. In other words, if $\varepsilon(D)=(D'_i)_{i\in\N}$, then $D'_i=D_i$ for all $i\in [m]$ and $D'_i(x_r)=0$ for all
 $i > m$ and all $r=1,\dots,d$. It is clear that $\varepsilon(\Id_A,\partial_s)$ coincides with the ``Taylor Hasse-Schmidt derivation''
 $\Delta^{(s)}$ defined in \ref{def:Taylor-HS}.
\end{nota}

Definition \ref{def:HS-integ} admits the following obvious logarihtmic version.

\begin{definicion} Let $I\subset A$ be an ideal and $n\in \overline{\N}$. We say that:
\begin{enumerate}
\item[1)] A $I$-logarithmic derivation $\delta\in \Der_k(\log I)$
 is {\em $I$-logarithmically $n$-integrable} if there is a $D\in
\HS_k(\log I;n)$ such that $D_1= \delta$.  A such $D$ will be called a {\em
$I$-logarithmic $n$-integral} of $\delta$. The set of $I$-logarithmic $k$-linear
derivations of $A$ which are $I$-logarithmically $n$-integrable will be denoted by
$\Ider_k(\log I;n)$. When $n=\infty$ it will be simply denoted by $\Ider_k(\log I)$.
\item[2)] A $I$-logarithmic Hasse--Schmidt derivation $D'\in\HS_k(\log I;m)$, with $m\leq n$,
is {\em $I$-logarithmically $n$-integrable} if there is a $D\in \HS_k(\log I;n)$
such that $\tau_{nm}D=D'$.  A such $D$ will be called a {\em $I$-logarithmic
$n$-integral} of $D'$. The set of $I$-logarithmically $n$-integrable $I$-logarithmic
Hasse-Schmidt derivations of $A$ over $k$ of length $m$  will be denoted by
$\IHS_k(\log I;m;n)$. When $n=\infty$ it will be simply denoted by $\IHS_k(\log
I;m)$.
\end{enumerate}
\end{definicion}

It is clear that the $\Ider_k(\log I;n)$ are $A$-submodules of $\Der_k(\log I)$ and
$\Der_k(\log I)=\Ider_k(\log I;1)\supset \Ider_k(\log I;2)\supset \cdots$
\begin{equation} \label{eq:intersection-iderlog}
\Ider_k(\log I) \subset \bigcap_{n\in\N_+} \Ider_k(\log I;n).
\end{equation}
It is also clear that the $\IHS_k(\log I;m;n)$ are subgroups of $\IHS_k(\log I;m)$,
stable by the $\sbullet$ operation, $\IHS_k(\log I;m) = \IHS_k(\log I;m;m) \supset
\IHS_k(\log I;m;m+1)\supset \cdots$ and
\begin{equation} \label{eq:intersection-iHSlog}
 \IHS_k(\log I;m) \subset \bigcap_{n\geq m} \IHS_k(\log I;m;n).
\end{equation}

The inclusions (\ref{eq:intersection-iderlog}) and (\ref{eq:intersection-iHSlog}) seem
not to be equalities in general (see question \ref{cuestion:1}). Nevertheless, we have the following proposition.

\begin{proposicion} \label{prop:inclusions-equalities} The following properties hold:
\begin{enumerate}
\item[1)] Let $n\geq 1$ be an integer. If any $k$-derivation of $A$ is $n$-integrable,
then any Hasse--Schmidt derivation $D\in\HS_k(A;m)$ is also $n$-integrable, for
all $m\leq n$.
\item[2)] If any $k$-derivation is $n$-integrable for all integers $n\geq 1$, then
any Hasse--Schmidt derivation $D\in\HS_k(A;m)$ is also $\infty$-integrable, for all
integers  $m\geq 1$.
\end{enumerate}
\end{proposicion}

\begin{prueba} For 1) we can mimic the proof of Proposition 1.4 in \cite{nar_2009}
by using Theorem 2.8 in \cite{magda_nar_hs} (see Remark 1.5 in \cite{nar_2009}). For
2), we apply 1) and we obtain a sequence $E^n\in\HS_k(A;n)$, $n\geq m$, with $E^m=D$
and $\tau_{n+1,n}E^{n+1}=E^n$ for all $n\geq m$. It is clear that the inverse limit
of the $E^n$ (see \ref{nume:truncation-inverse-limit}) is a $\infty$-integral of
$D$.
\end{prueba}

\begin{lema}  \label{lema:technical}
Assume that $R=k[x_1,\dots,x_d]$, $S\subset R$ is a multiplicative set and
$A=S^{-1}R$ or $A=k[[x_1,\dots,x_d]]$. Let $I\subset A$ be an ideal and $n\geq 1$ an
integer. Then, any Hasse--Schmidt derivation $D$ in the kernel of the group
homomorphism $\Pi_n$ (see \ref{num:HS-induced-by-HS-log}) is $I$-logarithmically
($\infty$-)integrable.
\end{lema}

\begin{prueba} Let us prove the proposition in the case $A=S^{-1}R$, the case
$A=k[[x_1,\dots,x_d]]$ being completely similar. Denote by
$\widetilde{\delta_r}:A\to A$ the induced derivation by the partial derivative
$\partial_r:R\to R$.
 We proceed by decreasing induction on $\ell(D)$ (see Definition \ref{def:ell}). If
$\ell(D)=n$, then $D$ is the identity and the result is clear. Let $m$ be an integer
with $0\leq m < n$ and suppose that any $D'\in\ker \Pi_n$ with $m+1\leq \ell(D')$ is
$I$-logarithmically integrable, and let $D\in\ker \Pi_n$ with $\ell(D)=m$, i.e. $D$
has the form $(\Id_A,0,\dots,0,D_{m+1},\dots,D_n)$ with $D_{m+1}\neq 0$, and so
$D_{m+1}$ must be a $k$-derivation. Since $D\in \ker \Pi_n$, we deduce that
$D_i(A)\subset I$ for all $i$. In particular, there are $a_1,\dots,a_d\in I$ such
that $D_{m+1}=\sum_{r=1}^d a_r \widetilde{\delta_r}$.

The $I$-logarithmic Hasse-Schmidt derivation $E=(a_1\sbullet
\widetilde{\Delta^{(1)}})\pcirc \cdots \pcirc (a_d\sbullet
\widetilde{\Delta^{(d)}})$ $\in \ker \Pi_{\infty}$ is an ($\infty$-)integral of
$D_{m+1}$. Let us consider $D'=D\pcirc (\tau_{\infty n} E[m+1])^{-1}\in \ker \Pi_n$.
It is clear that $\ell(D')\geq m+1$ and, by induction hypothesis, $D'$ is
$I$-logarithmically integrable. We conclude that $D=D'\pcirc (\tau_{\infty n}
E[m+1])$ is also $I$-logarithmically integrable.
\end{prueba}

\begin{nota} The proof of the above lemma shows that $\ker \Pi_n$ is
generated by the $n$-truncations of the $(a\sbullet E)[m]$, with $a\in I$,
$E\in\HS_k(A)$, $m\in [n]$. In fact, for $n=\infty$ we obtain that $\ker \Pi_\infty$
is the closure of subgroup of $\HS_k(\log I)$ generated by the $(a\sbullet E)[m]$,
with $a\in I$, $E\in\HS_k(A)$ and $m\in\N_+$, where we consider in $\HS_k(A)$ the
inverse limit topology of the discrete topologies in the $\HS_k(A;m)$, $m\in\N$ (see
\ref{nume:truncation-inverse-limit}). Namely, for $D\in\ker \Pi_\infty$, by the same
procedure as in the proof of the lemma we construct inductively a sequence
$E^q=(a_1^q\sbullet \widetilde{\Delta^{(1)}})\pcirc \cdots \pcirc (a_d^q\sbullet
\widetilde{\Delta^{(d)}})$, $q\geq 1$, $a_r^s\in I$, such that $\ell(D\pcirc
(F^q)^{-1})\geq q$, where $F^q=E^q[q]\pcirc \cdots \pcirc E^1[1]$. So $D\pcirc
(F^q)^{-1}$ tends to the identity element as $q \to \infty$ and $D$ is the limit of
$F^q$ as $q \to \infty$.
\end{nota}

\begin{proposicion} \label{prop:crit-integ-log-HS}
Assume that $R=k[x_1,\dots,x_d]$, $S\subset R$ is a multiplicative set and
$A=S^{-1}R$ or $A=k[[x_1,\dots,x_d]]$. Let $I\subset A$ be an ideal, $m\geq 1$ an
integer, $n\in\overline{\N}$ with $n\geq m$ and $E\in\HS_k(A/I;m)$. The following
properties are equivalent:
\begin{enumerate}
\item[(a)] $E$ is $n$-integrable.
\item[(b)] Any $D\in \HS_k(\log I;m)$ with
$\overline{D}=E$ is $I$-logarithmically
$n$-integrable.
\item[(c)] There is a $D\in \HS_k(\log I;m)$ with
$\overline{D}=E$ which is $I$-logarithmically $n$-integrable.
\end{enumerate}
\end{proposicion}

\begin{prueba} The implication (b) $\Rightarrow$ (c) is an obvious consequence of
Proposition \ref{prop:surjec-HS-log} and (c) $\Rightarrow$ (a) comes from
\ref{num:HS-induced-by-HS-log}. For the remaining implication (a) $\Rightarrow$ (b),
let $Z\in\HS_k(A/I;n)$ be an $n$-integral of $E$ and let $D\in\HS_k(\log I;m)$ be a
logarithmic Hasse-Schmidt derivation with $\overline{D}=E$. From Proposition
\ref{prop:surjec-HS-log}, there is a $U\in\HS_k(\log I;n)$ such that
$\overline{U}=Z$. Since $\overline{\tau_{nm}U} = \tau_{nm}\overline{U}= \tau_{nm}Z=
E = \overline{D}$, we have $D\pcirc (\tau_{nm}U)^{-1} \in \ker \Pi_m$ and so, by
Lemma \ref{lema:technical}, we deduce that $D$ is $I$-logarithmically
$n$-integrable.
\end{prueba}

\begin{corolario} \label{cor:crit-integ-log-HS} Under the hypotheses of Proposition \ref{prop:crit-integ-log-HS},
the map $\Pi_m:\IHS_k(\log I;m;n) \to \IHS_k(A/I;m;n)$ is surjective.
\end{corolario}

\begin{corolario} Under the hypotheses of Proposition \ref{prop:crit-integ-log-HS},
the following properties are equivalent:
\begin{enumerate}
\item[(a)] $\IHS_k(A/I;m;n)=\HS_k(A/I;m)$.
\item[(b)] $\IHS_k(\log I;m;n)=\HS_k(\log I;m)$.
\end{enumerate}
\end{corolario}

\begin{prueba} It is a straightforward consequence of the
proposition.
\end{prueba}

\begin{ejemplo} \label{ej:ncd} {\em (Normal crossings)} Let us take $\displaystyle f=\prod_{i=1}^e x_i\in A=k[x_1,\dots,x_d]$
and $I=(f)\subset A$. The $A$-module $\Ider_k(\log I)$ is generated by
$$\textstyle\left\{x_1\partial_1,\dots,x_e\partial_e,\partial_{e+1},\dots,\partial_d\right\}$$
and any of these $I$-logarithmic derivations are integrable $I$-logarithmically,
since $ \Delta^{(j)}, x_i\sbullet \Delta^{(i)} \in \HS_k(\log I)$ for $i=1,\dots,e$
and $j=e+1,\dots,n$. In particular $\Ider_k(\log I)=\Der_k(\log I)$ and
$\Ider_k(A/I)=\Der_k(A/I)$.
\end{ejemplo}

\begin{proposicion} \label{prop:algorit-1}  Let $A$ be an arbitrary $k$-algebra,
$I\subset A$ an ideal with generators $f_l$, $l\in L$, and $n\geq 1$ an integer. Let $D\in\HS_k(\log I;n)$ be a $I$-logarithmic Hasse--Schmidt derivation and assume that $D$ is $(n+1)$-integrable and let $(\Id_A,D_1,\dots,D_n,D_{n+1})\in\HS_k(A;n+1)$ be an $(n+1)$-integral of $D$. The following properties are
equivalent:
\begin{enumerate}
\item[(a)] $D$ is $I$-logarithmically $(n+1)$-integrable.
\item[(b)] There is a derivation $\delta\in\Der_k(A)$ such that $D_{n+1}(f_l) + \delta(f_l) \in I$ for all $l\in L$
\end{enumerate}
\end{proposicion}

\begin{prueba} It comes from the fact that any other $(n+1)$-integral of $D$ must be of the form $(\Id_A,D_1,\dots,D_n,D_{n+1}+\delta)$ with $\delta\in \Der_k(A)$.
\end{prueba}

\begin{corolario} \label{cor:algorit-2} Assume
that $A=k[x_1,\dots,x_d]$ or $A=k[[x_1,\dots,x_d]]$. Let $I=(f_1,\dots,f_p)\subset A$ be an ideal and $n\geq 1$ an integer. Let $D\in\HS_k(\log I;n)$ be a $I$-logarithmic Hasse--Schmidt derivation and let us consider its integral $D'=\varepsilon(D)$ (see remark \ref{nota:taylor-HS}). The following properties are
equivalent:
\begin{enumerate}
\item[(a)] $D$ is $I$-logarithmically $(n+1)$-integrable.
\item[(b)] There are $\alpha_r, a_{st}\in A$, $r=1,\dots,d$, $s,t=1,\dots,p$, such that
$$
D'_{n+1}(f_s) = \alpha_1 \left(f_s\right)'_{x_1} + \cdots + \alpha_d \left(f_s\right)'_{x_d} + a_{s1}f_1+\cdots+a_{sp}f_p\quad \forall s=1,\dots,p.
$$
\end{enumerate}
Moreover, if (b) holds, an explicit $I$-logarithmic $(n+1)$-integral of $D$ is given by $(\Id_A,D_1,\dots,D_n,D'_{n+1}-\delta)$, with $\delta=\sum_{r=1}^d \alpha_r
\partial_r$.
\end{corolario}

\begin{nota} \label{nota:comput-1} (1) In the case of a ``computable'' base ring $k$ (for instance, any finitely
generated extension of $\ZZ,\QQ$ or of any finite field) and a finitely presented
$k$-algebra $A$, Proposition \ref{prop:crit-integ-log-HS} and Corollary
\ref{cor:algorit-2} give an effective way to decide whether a given Hasse--Schmidt
derivation $D\in\HS_k(A;n)$ of finite length $n$ is $(n+1)$-integrable or not and,
if yes, to compute an explicit $(n+1)$-integral of $D$.
\smallskip

\noindent (2) Nevertheless, the question of deciding whether a given Hasse--Schmidt
derivation $D\in\HS_k(A;n)$ of finite length $n$ is $(n+r)$-integrable or not, with
$r\geq 2$, is much more involved. First of all, we cannot proceed ``step by step'',
since $D$ can be $(n+r)$-integrable and simultaneously admit an $(n+1)$-integral
which is not $(n+r)$-integrable (cf. example 3.7 in \cite{nar_2009}). On the other
hand, the condition of $(n+r)$-integrability of $D$, $r\geq 2$, gives rise to
nonlinear equations which seem not obvious to treat in general with the currently
available methods, either theoretical or computational (see for instance Lemmas
\ref{lema:ej-nolineal}, \ref{lema:ej-nolineal-3}, \ref{lema:ej-nolineal-4}).
\smallskip

\noindent (3)  The following example is a very particular case of a general result, but it also serves to
illustrate the nonlinear nature of integrability and the difficulties that come
from: Let $A=k[x_1,\dots,x_d]$, $f\in A$, $I=(f)$ and
$\delta=\sum_{r=1}^d a_r
\partial_r$ any $k$-derivation of $A$. The following properties
are equivalent:
\begin{enumerate}
\item[(a)] $\delta$ is a $I$-logarithmic derivation and it is $I$-logarithmically
$2$-integrable.
\item[(b)] $\sum_{r=1}^d  f'_{x_r} a_r \in I$ and
$\sum_{|\alpha|=2} \Delta^{(\alpha)}(f)\ \underline{a}^\alpha \in
(f,f'_{x_1},\dots,f'_{x_d})$.
\end{enumerate}
So, in order to compute a system of generators of the $A$-module $\Ider_k(\log
I;2)$, we have to deal with nonlinear homogeneous equations of degree $2$ (see examples in sections \ref{subsect:cusp-2-3}, \ref{subsect:cusp2-Z}).
\end{nota}

\subsection{Jacobians and integrability} \label{subsect:jacob}

Let $k$ be an arbitrary (commutative) ring and assume that $R=k[x_1,\dots,x_d] $ or  $R=k[[x_1,\dots,x_d]] $. Let  $I=(f_1,\dots,f_u)\subset R$ be a finitely generated ideal and $A=R/I$.
For each $e=1,\dots,\min\{d,u\}$ let $J^0_e$ be the ideal generated by  all the $e\times e$ minors of the Jacobian matrix $(\partial f_j/\partial x_i)$, and $J_e = (J^0_e+I)/I$.
 We have $J_1 \supset J_2 \supset \cdots$. Let $c$ be the maximum index $e$ with $J_e\neq 0$
 (or equivalently with $J^0_e\nsubseteq I$), in case it exists. The ideal $J_c$ only depends
 on the $k$-algebra $A$ and is called the {\it Jacobian ideal} of $A$ over $k$ and denoted by
 $J_{A/k}$. It is nothing else but the smallest non-zero Fitting ideal of the module of
 $k$-differentials $\Omega_{A/k}$ (see \cite{lipman_1969}).

\begin{proposicion} \label{prop:Jder-in-Ider}  Under the above hypotheses, any $\delta\in \Der_k(\log I) \cap (J^0_c+I) \Der_k(R)$ is $I$-logarithmically integrable.
\end{proposicion}

\begin{prueba} The proof follows the same lines that the proof of Theorem 11 in \cite{mat-intder-I}. Let us write $J^0=J^0_c$. Since $I \Der_k(R) \subset \Ider_k(\log I)$, we can assume that
$\delta=\sum_{r=1}^d c_{r1} \partial_r$ with $c_{r1}\in J^0$.
Let us consider $D^1=(\Id_A,\delta)\in\HS_k(\log I;1)$ and
$E^1=\varepsilon(D^1)\in\HS_k(R;\infty)$ (see \ref{nota:taylor-HS}). We have that
$E^1_2 = \sum_{|\alpha|=2} \left(\prod_{r=1}^d c_{r1}^{\alpha_r}\right)
\Delta^{(\alpha)}\in (J^0)^2\diff_{R/k}$, and so $E^1_2(f_j)\in (J^0)^2$ for all
$j=1,\dots,u$. From Lemma \ref{lema:gene-jacob-matsu} there is
$(c_{12},\dots,c_{d2})\in R^d$, with $c_{r2}\in J^0$, such that
$$ (c_{12},\dots,c_{d2}) \left((\partial f_j/\partial x_i)_{i,j}\right) \equiv (E^1_2(f_1),
\dots,E^1_2(f_u)) \mod I,$$ i.e. $E^1_2(f_j)-\sum_{r=1}^d c_{r2}(f_j)'_{x_r}\in I$,
and so we deduce that $D^1$ is $I$-logarithmically $2$-integrable, an
$I$-logarithmic $2$-integral being $D^2=(\Id_A,\delta,D^2_2)$ with $D^2_2=E^1_2 -
\sum_{r=1}^d c_{r2} \partial_r\in J^0 \diff_{R/k}$ (see Corollary
\ref{cor:algorit-2}). \smallskip

Assume that we have found a $D^m=(\Id_A,\delta,D^2_2,\dots,D^m_m)\in \HS_k(\log
I;m)$ with $D^s_s \in J^0 \diff_{R/k}$, $s=1,\dots,m$, hence with
$c_{rs}:=D^s_s(x_r)\in J^0$, $r=1,\dots,d$. Let us consider $E^m=\varepsilon(D^m)\in
\HS_k(R;\infty)$. From Proposition \ref{prop:formulon}, 2) we deduce that 
$E^m_{m+1} \in (J^0)^2 \diff_{A/k}$ and so $E^m_{m+1}(f_j)\in (J^0)^2$ for all
$j=1,\dots,u$. From Lemma \ref{lema:gene-jacob-matsu}, there is
$(c_{1,m+1},\dots,c_{d,m+1})\in R^d$, with $c_{r,m+1}\in J^0$, such that
$$ (c_{1,m+1},\dots,c_{d,m+1}) \left((\partial f_j/\partial x_i)_{i,j}\right) \equiv (E^m_{m+1}(f_1),\dots,E^m_{m+1}(f_u)) \mod I,
$$
i.e. $E^m_{m+1}(f_j)- \sum_{r=1}^d c_{r,m+1}(f_j)'_{x_r} \in I$, and so we deduce
again that $D^m$ is $I$-logarithmically $m+1$-integrable, an $I$-logarithmic
$(m+1)$-integral being $D^{m+1}=(\Id_A,\delta,D^2_2,\dots,\D^m_m,D^{m+1}_{m+1})$
with $D^{m+1}_{m+1}=E^m_{m+1} - \sum_{r=1}^d c_{r,m+1} \partial_r\in J^0
\diff_{R/k}$ (see Corollary \ref{cor:algorit-2}). \smallskip

In that way, we construct inductively the $D^m_m$, $m\geq 2$, such that $(\Id_A,\delta,D^2_2,\dots)\in \HS_k(\log I;\infty)$ and so $\delta$ is $I$-logarithmically integrable.
\end{prueba}

\begin{lema} \label{lema:gene-jacob-matsu} Let $\mathbf{X}=(X_{ij})$, $i=1,\dots,d$, $j=1,\dots,u$, be variables, $W=\ZZ[\mathbf{X}]$, $\mathfrak{a}_{e}\subset W$ the ideal generated by the $e\times e$ minors of $\mathbf{X}$ and $U=W/\mathfrak{a}_{c+1}$. Then, for each $c\times c$ minor $\mu$ of $\mathbf{X}$ and for each $j=1,\dots,u$, the system
$$     (u_1,\dots,u_d) \mathbf{X} = (0,\dots,0,\stackrel{\scriptsize\underbrace{j}}{\mu},0,\dots,0)
$$
has a solution in $U$.
\end{lema}

\begin{prueba} We know that $U$ is an integral domain (cf. \cite{bruns-vetter-1988}, Theorem (2.10) and Remark (2.12)).
Denote by $K$ its field of fractions and by $\pi:W\to U$ the natural projection. The
lemma is an easy consequence of the fact that the matrix $\pi(\mathbf{X})\otimes K$
has rank $c$.
\end{prueba}

The following corollary of Proposition \ref{prop:Jder-in-Ider} generalizes Theorem 11 in \cite{mat-intder-I}, which was only stated and proved for $k$ a perfect field.

\begin{corolario} Under the above hypotheses, we have
$$J_{A/k}\subset \ann_A\left(\Der_k(A)/\Ider_k(A)\right).$$
\end{corolario}

The proof of the following result is similar to the proof of Proposition \ref{prop:Jder-in-Ider}.

\begin{proposicion}
Let $f\in R$, $I=(f)$, and $J^0=(f'_{x_1},\dots,f'_{x_d})$ the gradient ideal.
If $\delta:R\to R$ is a $I$-logarithmic $k$-derivation with $\delta\in J^0\Der_k(R)$, then $\delta$ admits a $I$-logarithmic integral $D\in\HS_k(\log I)$
 with $D_i(f)=0$  for all $i>1$. In particular, if $\delta(f)=0$, the integral $D$ can be taken with $\Phi_D(f)=f$.
\end{proposicion}

\noindent \numero \label{nume:traves} We quote here Theorem 1.2 in
\cite{traves-2000}: Let $I\subset A=k[x_1,\dots,x_d]$ be an ideal generated by
quasi-homogeneous polynomials with respect to the weights $w(x_r)\geq 0$. Then, the
Euler vector field $\chi = \sum_{r=0}^d w(x_r) \partial_r$ is $I$-logarithmically
($\infty$-)integrable. In fact, a $I$-logarithmic integral of $\chi$ is the
Hasse--Schmidt derivation associated with the map $A \to A[[t]]$ given by
\begin{eqnarray*}
&x_r  \mapsto x_r\left(\frac{1}{1-t}\right)^{w(x_r)},\quad r=1,\dots,d.&
\end{eqnarray*}

\begin{proposicion} \label{prop:isol-sing-general} Let $f\in A=k[x_1,\dots,x_d]$ be a quasi-homogeneous
polynomial with respect to the weights $w(x_r)> 0$ and $I=(f)\subset A$. Assume that
the weight of $f$ is a unit in $k$ and that all the partial derivatives of $f$ are
non-zero and form a regular sequence. Then $\Der_k(\log I)=\Ider_k(\log I)$.
\end{proposicion}

\begin{prueba} From the hypotheses we deduce that the $A$-module
$\Der_k(\log I)$ is generated by the Euler vector field $\chi$ and the crossed
derivations $\theta_{rs}=f'_{x_s}\partial_r - f'_{x_r}\partial_s$, $1\leq r < s\leq
d$. But $\chi$ is $I$-logarithmically integrable by \ref{nume:traves} and
$\theta_{rs}$ is $I$-logarithmically integrable by Proposition
\ref{prop:Jder-in-Ider}.
\end{prueba}

\subsection{Behaviour of integrability under localization}

Throughout this section, $k$ will be an arbitrary commutative ring.
\medskip

The proof of the following proposition is clear from \ref{nume:basic-localiz}.

\begin{proposicion} \label{prop:coh-1} Let $A$ be a $k$-algebra, $S\subset A$ a multiplicative
set, $\mathfrak{a}\subset A$ be an ideal, $m\geq 1$ an integer, $n\in\overline{\N}$
with $n\geq m$ and $D\in\HS_k(\log \mathfrak{a};m)$. If $D$
$\mathfrak{a}$-logarithmically $n$-integrable, then
$\widetilde{D}\in\HS_k(S^{-1}A;m)$ is $(S^{-1}\mathfrak{a})$-logarithmically
$m$-integrable. In particular, the map $\Theta_m$ sends $\IHS_k(\log
\mathfrak{a};m;n)$ to $\IHS_k(\log (S^{-1}\mathfrak{a});m;n)$.
\end{proposicion}

The two following propositions are straightforward consequences of Proposition
\ref{prop:surjec-HS-log-localiz-poly} and Corollary
\ref{cor:surjec-HS-log-localiz-fp} respectively.

\begin{proposicion} \label{prop:surjec-IHS-log-localiz-poly}
Assume that $A=k[x_1,\dots,x_d]$ and let $S\subset A$ be a multiplicative set and
$\mathfrak{a}=(f_1,\dots,f_u)\subset A$ be a finitely generated ideal. Then, for any
integers $m\geq q\geq 1$, the map
$$ (s,F)\in S \times \IHS_k(\log \mathfrak{a};q;m)\mapsto \frac{1}{s}\sbullet
\Theta_q(F)
\in\IHS_k(\log (S^{-1}\mathfrak{a});q;m)$$ is surjective.
\end{proposicion}

\begin{proposicion} \label{prop:surjec-IHS-localiz-fp} Assume that $A$ is a finitely presented $k$-algebra
and let $T\subset A$ be a multiplicative set. Then, for any integers $m\geq q\geq 1$
the map
$$ (t,G)\in T \times \IHS_k(A;q;m)\mapsto \frac{1}{t}\sbullet \Theta_q(G)
\in\IHS_k(T^{-1}A;q;m)$$ is surjective.
\end{proposicion}

Proposition \ref{prop:surjec-IHS-localiz-fp} can be also obtained form Proposition
\ref{prop:surjec-IHS-log-localiz-poly} and Corollary \ref{cor:crit-integ-log-HS}.

\begin{corolario} \label{cor:coh-2}
Assume that $A=k[x_1,\dots,x_d]$ and let $S\subset A$ be a multiplicative set,
 $\mathfrak{a}=(f_1,\dots,f_u)\subset A$ be a
finitely generated ideal. Then, for any integer $m\geq 1$ the canonical map
$$ \frac{\delta}{s}\in S^{-1} \Ider_k(\log \mathfrak{a};m) \mapsto
\frac{1}{s}\widetilde{\delta}\in \Ider_k(\log
(S^{-1}\mathfrak{a});m)$$ is an isomorphism of $(S^{-1}A)$-modules.
\end{corolario}

\begin{prueba} The injectivity is a consequence of the fact that, under the
above assumptions, the canonical map $S^{-1}\Der_k(A) \to \Der_k(S^{-1}A)$ is an
isomorphism. The surjectivity is given by Proposition
\ref{prop:surjec-IHS-log-localiz-poly} in the case $q=1$.
\end{prueba}

\begin{corolario} \label{cor:ider-localiz_fp} Assume that $A$ is a finitely presented $k$-algebra and let $T\subset A$
be a multiplicative set. Then, for any integer $m\geq 1$ the canonical map
$$  T^{-1} \Ider_k(A;m) \to \Ider_k(T^{-1}A;m)$$
is an isomorphism of $(T^{-1}A)$-modules.
\end{corolario}

\begin{prueba} The injectivity goes as in the proof
of Corollary \ref{cor:coh-2}. The surjectivity is given by Proposition
\ref{prop:surjec-IHS-localiz-fp} in the case $q=1$.
\end{prueba}

\begin{teorema} \label{teo:criter-loc-ider-fp}
Assume that $A$ is a finitely presented $k$-algebra, $m\geq 1$ is an integer and let
$\delta\in \Der_k(A)$. The following properties are equivalent:
\begin{enumerate}
\item[(a)] $\delta\in \Ider_k(A;m)$.
\item[(b)] $\delta_{\mathfrak{p}}\in \Ider_k(A_{\mathfrak{p}};m)$ for all $\mathfrak{p}\in\Spec A$.
\item[(c)] $\delta_{\mathfrak{m}}\in \Ider_k(A_{\mathfrak{m}};m)$ for all $\mathfrak{m}\in\Specmax A$.
\end{enumerate}
\end{teorema}

\begin{prueba} The implication (a) $\Rightarrow$ (b) is a consequence of Proposition
\ref{prop:coh-1}. The implication (b) $\Rightarrow$ (c) is obvious. For the
remaining implication (c) $\Rightarrow$ (a), assume that property (c) holds. Then,
by Corollary \ref{cor:ider-localiz_fp}, for each $\mathfrak{m}\in\Specmax A$ there is a
$f^\mathfrak{m} \in A - \mathfrak{m}$ and a $\zeta^\mathfrak{m} \in \Ider_k(A;m)$
such that $f^\mathfrak{m}\delta_{\mathfrak{m}} =
\left(\zeta^\mathfrak{m}\right)_\mathfrak{m}$, and so there is a $g^\mathfrak{m} \in A -
\mathfrak{m}$ such that $g^\mathfrak{m}f^\mathfrak{m}\delta =
g^\mathfrak{m}\zeta^\mathfrak{m}$. Since the ideal generated by the $g^\mathfrak{m}
f^\mathfrak{m}$, $\mathfrak{m}\in\Specmax A$,  must be the total ideal, we deduce
the existence of a finite number of $\mathfrak{m}_i\in\Specmax A$ and $a_i\in A$,
$1\leq i\leq n$, such that $1=a_1 g_1f_1 +\cdots + a_n g_nf_n$, with $f_i=
f^{\mathfrak{m}_i}$, $g_i= g^{\mathfrak{m}_i}$, and so
$$ \delta = \sum_{i=1}^n a_i g_i f_i \delta = \sum_{i=1}^n a_i
g_i\zeta^{\mathfrak{m}_i}$$ is $m$-integrable.
\end{prueba}

\begin{corolario} \label{cor:Ider_fp_schemes}
Let $f:X \to S$ be a locally finitely presented morphism of schemes. For each
integer $n\geq 1$ there is a quasi-coherent sub-sheaf $\fIder_S(\OO_X;n)\subset
\fDer_S(\OO_X)$ such that, for any affine open sets $U=\Spec A \subset X$ and
$V=\Spec k \subset S$, with $f(U)\subset V$, we have $\Gamma(U,\fIder_S(\OO_X;n))
=\Ider_k(A;n)$ and $\fIder_S(\OO_X;n)_p = \Ider_{\OO_{S,f(p)}}(\OO_{X,p};n)$ for
each $p\in X$. Moreover, if $S$ is locally noetherian, then $\fIder_S(\OO_X;n)$ is a
coherent sheaf.
\end{corolario}

\begin{prueba} For each open set $U\subset X$, we define
$$\Gamma(U,\fIder_S(\OO_X;n)) =\{\delta\in \Gamma(U,\fDer_S(\OO_X))\ |\ \delta_p \in  \Ider_{\OO_{S,f(p)}}(\OO_{X,p};n)\ \forall p\in U\}.$$
The behaviour of $\fIder_S(\OO_X;n)$ on affine open sets and its quasi-coherence is
a straightforward consequence of Theorem \ref{teo:criter-loc-ider-fp}.
\end{prueba}

\subsection{Testing the integrability of derivations}

In this section $k$ will be an arbitrary commutative ring and $A$ an arbitrary $k$-algebra.

\begin{definicion} Let $n\geq m > 1$ be integers and $D\in\HS_k(A;n)$.
We say that $D$ is {\em $m$-sparse} if $D_i=0$ whenever $i\notin\N m$. We say that
$D$ is {\em weakly $m$-sparse} if $\tau_{n,qm}D$ is $m$-sparse, where
$q=\left\lfloor\frac{n}{m}\right\rfloor$. The set of $m$-sparse (resp. weakly
$m$-sparse) Hasse--Schmidt derivations in $\HS_k(A;n)$ will be denoted by
$\HS_k^{m-sp}(A;n)$ (res. $\HS_k^{m-wsp}(A;n)$).
\end{definicion}

The proof of the following proposition is easy and its proof is left up to the reader.

\begin{proposicion} \label{prop:sparse} Let $n\geq m > 1$ be integers, $q=\left\lfloor\frac{n}{m}\right\rfloor$ and $r=n-qm$. The following properties hold:
\begin{enumerate}
\item[1)] $\HS_k^{m-sp}(A;n)$ and $\HS_k^{m-wsp}(A;n)$ are subgroups of $\HS_k(A;n)$.
\item[2)] For any $D\in\HS_k(A;q)$ and any $\underline{\delta}=(\delta_1,\dots,\delta_r)\in \Der_k(A)^r$, the sequence
$$ \Theta(D,\underline{\delta})=(\Id_A,0,\dots,0,\stackrel{\scriptsize\underbrace{m}}{D_1},0,\dots,0,\stackrel{\scriptsize\underbrace{2m}}{D_2},0,\dots,0,\stackrel{\scriptsize\underbrace{qm}}{D_q},
\stackrel{\scriptsize\underbrace{qm+1}}{\delta_1},\dots,\stackrel{\scriptsize\underbrace{n}}{\delta_r})$$
is a weakly $m$-sparse Hasse--Schmidt derivation of $A$ (over $k$) of length $n$ and the map
$\Theta: \HS_k(A;q) \times \Der_k(A)^r \to \HS_k^{m-wsp}(A;n)$
is an isomorphism of groups.
\end{enumerate}
\end{proposicion}

\begin{teorema} \label{th:main} Let $n\geq 1$ be an integer.
The following assertions hold:
\begin{enumerate}
\item[1)] If $n$ is odd and $\Ider_k(A;q) = \Der_k(A)$, with $q= \frac{n+1}{2}$,
then any $D\in\HS_k(A;n)$ with $D_1=0$ is $(n+1)$-integrable.
\item[2)] If $n$ is even and $\Ider_k(A;p) = \Der_k(A)$, with $p=\left\lfloor
\frac{n+1}{3}\right\rfloor$, then any $D\in\HS_k(A;n)$ with $D_1=0$ is
$(n+1)$-integrable.
\end{enumerate}
\end{teorema}

\begin{prueba} 1) Since $D_1=0$ we have $1\leq\ell(D)\leq n$.
If $n=1$, then $D$ is the identity and the result is clear. Assume $n\geq 3$ and so
$q\geq 2$. Let us proceed by decreasing induction on $\ell(D)$. If $\ell(D) =n$ then
$D$ is the identity and the result is clear. Let $m$ be an integer with $1\leq m <
n$ and suppose that any $D'\in\HS_k(A;n)$ with $m+1\leq \ell(D')$ is
$(n+1)$-integrable. Let $D\in\HS_k(A;n)$ be a Hasse--Schmidt derivation with $\ell(D)=m$, i.e.
$$ D=(\Id_A,0,\dots,0,D_{m+1},\dots,D_n)\quad\text{with}\ D_{m+1}\neq 0.$$
Since $\tau_{n,m+1}D$ is $(m+1)$-sparse, we can apply Proposition \ref{prop:sparse}, 2) and deduce that
$D_{m+1}$ is a derivation and so, by hypothesis, it must be
$q$-integrable. Let $E\in\HS_k(A;q)$ be a $q$-integral of $D_{m+1}$. We have that
$q(m+1) \geq 2q =n+1$ and so $F=\tau_{q(m+1),n}(E[m+1])$ is $(n+1)$-integrable, an $(n+1)$-integral being
$\tau_{q(m+1),n+1}(E[m+1])$, and
has the form
$$ F=(\Id_A,0,\dots,0,\stackrel{\scriptsize\underbrace{m+1}}{D_{m+1}},0,\dots,F_n).
$$
It is clear that for $D' = F^{-1}\pcirc D$ we have $D'_1=\dots=D'_{m+1}=0$, and so
$\ell(D')\geq m+1$. The induction hypothesis implies that $D'$ is $(n+1)$-integrable
and we conclude that $D=F\pcirc D'$ is also $(n+1)$-integrable. \smallskip

\noindent 2) If $n=2$, then $D=(\Id_A,0,D_2)$ and obviously $(\Id_A,0,D_2,0)$ is a
$3$-integral of $D$. Assume that $n$ is even $\geq 4$, and let us write $n=2q$,
$q\geq 2$, and $n+1=3p+r$ with $0\leq r < 3$, $p\geq 1$. Since $\tau_{n3}D$ is weakly $2$-sparse, we deduce that
$D_3$ must be a derivation (see Proposition \ref{prop:sparse}) and so, by hypothesis, it is $p$-integrable. Let
$E^3\in\HS_k(A;p)$ be a $p$-integral of $D_3$. It is clear that (see Proposition \ref{prop:sparse})
$$ F^3=(\Id_A,0,0,\stackrel{\scriptsize\underbrace{3}}{E^3_1},0,0,\stackrel{\scriptsize\underbrace{6}}{E^3_2},0,\dots,0,
\stackrel{\scriptsize\underbrace{3p}}{E^3_p},0,0)$$ is a $(3p+2)$-integral of
$E^3[3]$, and since $3p+2\geq n+1$, $G^3=\tau_{3p+2,n}F^3$ is $(n+1)$-integrable and
$(G^3)^{-1}\pcirc D$ has the form $(\Id_A,0,D_2,0,\dots)$.

Assume that we have found $G^3,G^5,\dots,G^{2s-1}\in\HS_k(A;n)$, all of them
$(n+1)$-integrable, with $3\leq 2s-1 < n$, such that $(G^{2s-1})^{-1}\pcirc \cdots
\pcirc (G^3)^{-1}\pcirc D$ has the form
$$ D'=(\Id_A,0,\stackrel{\scriptsize\underbrace{2}}{D'_2},0,\stackrel{\scriptsize\underbrace{4}}{D'_4},0,\dots,0,
\stackrel{\scriptsize\underbrace{2s}}{D'_{2s}},D'_{2s+1},\dots,D'_n).$$ If $2s=n$,
we already have what we are looking for. If $2s<n$, then $D'_{2s+1}$
is a derivation (see Proposition \ref{prop:sparse}) and so, by hypothesis, it is $p$-integrable. Let
$E^{2s+1}\in\HS_k(A;p)$ be a $p$-integral of $D'_{2s+1}$. Let us consider
$F^{2s+1}=E^{2s+1}[2s+1]\in\HS_k(A;p(2s+1))$. Since $p(2s+1)\geq 5p \geq 3p+2\geq
n+1$, $G^{2s+1}:= \tau_{p(2s+1),n}F^{2s+1}$ is $(n+1)$-integrable and
$(G^{2s+1})^{-1}\pcirc D'$ has the form
$$ D''=(\Id_A,0,\stackrel{\scriptsize\underbrace{2}}{D''_2},0,\stackrel{\scriptsize\underbrace{4}}{D''_4},0,\dots,0,
\stackrel{\scriptsize\underbrace{2s}}{D''_{2s}},0,\stackrel{\scriptsize\underbrace{2s+2}}{D''_{2s+2}},\dots,D''_n).$$

We conclude with the existence of $G^3,G^5,\dots,G^{n-1}\in\HS_k(A;n)$, all of them
$(n+1)$-integrable, such that $ H = (G^{n-1})^{-1}\pcirc G^{n-3} \cdots \pcirc
(G^3)^{-1}\pcirc D\in\HS_k(A;n)$ ($n=2q$) is $2$-sparse. From Proposition \ref{prop:sparse} again we deduce that
$H$ is $(n+1)$-integrable, and so $D$ is also $(n+1)$-integrable.
\end{prueba}

\begin{definicion} For each integer $n\geq 1$, let us define
$$ \rho(n)= \left\{\begin{array}{ll}  \frac{n+1}{2} & \text{if $n$ is
odd}\\ \left\lfloor \frac{n+1}{3} \right\rfloor & \text{if $n$ is even.}\end{array}
\right.
$$
\end{definicion}
Notice that $\rho(n) < n$ for all $n\geq 2$.

\begin{corolario} \label{cor:main} Let $n\geq 1$ be an integer, and assume that
$\Ider_k(A;\rho(n)) = \Der_k(A)$. Then, for any $n$-integrable derivation
$\delta\in\Ider_k(A;n)$, the following properties are equivalent:
\begin{enumerate}
\item[(a)] Any $n$-integral of $\delta$ is
$(n+1)$-integrable.
\item[(b)] There is an $n$-integral of $\delta$ which is
$(n+1)$-integrable.
\end{enumerate}
\end{corolario}

\begin{prueba} Assume that $E\in\HS_k(A;n+1)$ is an $(n+1)$-integral of $\delta$ and let $D\in\HS_k(A;n)$ be any $n$-integral of $\delta$. The $1$-component of
$F=D\pcirc(\tau_{n+1,n}E)^{-1}$ vanishes and so, by Theorem \ref{th:main}, $F$ is $(n+1)$-integrable. We deduce that
$D=F\pcirc \tau_{n+1,n}E$ is also $(n+1)$-integrable.
\end{prueba}

\subsection{Algorithms} \label{subsect:algo}

Let $k$ be a ``computable'' base ring $k$ (for instance, any finitely generated
extension of $\ZZ,\QQ$ or of any finite field), $f_1,\dots,f_p\in
A=k[x_1,\dots,x_d]$ and $I=(f_1,\dots,f_p)$. The starting point is the computation
of a system of generators $\{\delta^1,\dots,\delta^q\}$ of $\Der_k(\log I)$.
\medskip

The following algorithm decides whether the equality
$$\Der_k(\log I)\stackrel{\text{?}}{=}\Ider_k(\log I;2)\quad \left(\Leftrightarrow \Der_k(A/I)\stackrel{\text{?}}{=}\Ider_k(A/I;2)\right)$$
 is true or not, and if yes, returns a $2$-integral for each generator of $\Der_k(\log I)$.
\bigskip

\noindent {\sc ALGORITHM--1:}
\begin{description}
\item[Step 1:] For each $j=1,\dots,q$, apply Corollary \ref{cor:algorit-2} as explained in remark \ref{nota:comput-1}, (1) to decide whether $\delta^j$ is $I$-logaritmically $2$-integrable or not, and if yes to compute a $I$-logarithmic $2$-integral $D^{j,2}$ of $\delta^j$.
\item[Step 2:] (Y) If the answer in Step 1 is YES for all $j=1,\dots,q$, then save the $I$-logarithmic $2$-integrals $D^{1,2},\dots,D^{q,2}$ and answer ``THE EQUALITY $\Der_k(\log I)=\Ider_k(\log I;2)$ IS TRUE''.\\
    (N) If the answer in step 1 is NOT for some $j=1,\dots,q$, then answer ``THE EQUALITY $\Der_k(\log I)=\Ider_k(\log I;2)$ IS FALSE''.
\end{description}
\medskip

Assume that we have an {\sc ALGORITHM--(N-1)} to decide whether the equality
$$\Der_k(\log I)\stackrel{\text{?}}{=}\Ider_k(\log I;N)\quad \left(\Leftrightarrow \Der_k(A/I)\stackrel{\text{?}}{=}\Ider_k)A/I;N)\right)$$
is true or not, and if yes, to compute an $N$-integral for each generator of $\Der_k(\log I)$.
\medskip

\noindent {\sc ALGORITHM--N:}
\begin{description}
\item[Step 1:] Apply {\sc ALGORITHM--(N-1)}, and if the answer is NOT, then STOP and answer ``THE EQUALITY $\Der_k(\log I)=\Ider_k(\log I;N+1)$ IS FALSE''.\\
If the answer to {\sc ALGORITHM--(N-1)} is YES, keep the computed $I$-logarithmic $N$-integrals $D^{1,N},\dots,D^{q,N}$ of $\delta^1,\dots,\delta^q$  and go to step 2.
\item[Step 2:] For each $j=1,\dots,q$, apply Corollary \ref{cor:algorit-2} as explained in remark \ref{nota:comput-1}, (1) to decide whether $D^{j,N}$ is $I$-logaritmically $(N+1)$-integrable or not, and if yes to compute a $I$-logarithmic $(N+1)$-integral $D^{j,N+1}$ of $D^{j,N}$.
\item[Step 3:] (Y) If the answer in Step 2 is YES for all $j=1,\dots,q$, then save the $I$-logarithmic $(N+1)$-integrals $D^{1,N+1},\dots,D^{q,N+1}$ and answer ``THE EQUALITY $\Der_k(\log I)=\Ider_k(\log I;N+1)$ IS TRUE''.\\
    (N) If the answer in Step 2 is NOT for some $j=1,\dots,q$, then answer ``THE EQUALITY $\Der_k(\log I)=\Ider_k(\log I;N+1)$ IS FALSE''.
\end{description}
\medskip

Corollary \ref{cor:main} is the key point for the correctness of Step 3, (N).

\section{Examples and questions}

We have used Macaulay 2 \cite{M2} for the preliminary computations needed in the following examples.

\subsection{The cusp $x^2+y^3$ in characteristic $2$ or $3$} \label{subsect:cusp-2-3}

Let $k$ be a base ring containing the field $\mathbb{F}_p$, $p>0$, and $f=x^2+y^3 \in R=k[x,y]$. Let $I=(f)$ and $A=k[x,y]/I$.
The computation of $\Ider_k(A;\infty)$ has been treated in \cite{mat-intder-I}, example 5. Here we are interested in the computation of $\Ider_k(A;m)$, $m\geq 2$.
\medskip

Let start with $p=2$. Then the Jacobian ideal of $f$ is
$J=(y^2,f)=(x^2,y^2)$.
\medskip

The module $\Der_k(\log I)$ is free with basis $\{\dpar{x}, f\dpar{y}\}$. It is
clear that $f\dpar{y}$ is $I$-logarithmically ($\infty$-)integrable. Let $g\in R$ be
a polynomial. From Corollary \ref{cor:algorit-2}, we have that $g\dpar{x}$ is
$I$-logarithmically $2$-integrable if and only if $g^2\in J$. Since $\{g\in R\ |\
g^2 \in J\} = (x,y)$, we deduce that $\{x\dpar{x},y\dpar{x},f\dpar{y}\}$ is a system
of generators of $\Ider_k(\log I;2)$.
\medskip

The derivation $x\dpar{x}$ is the Euler vector field for the weights $w(x)=3,
w(y)=2$. From \ref{nume:traves} we know that $x\dpar{x}$ is $I$-logarithmically
($\infty$-)integrable. \medskip

 Let $c\in R$ be an arbitrary polynomial and $\delta=cy\dpar{x}$. A $I$-logarithmic $2$-integral of $\delta$ is determined by the $k$-algebra map
$$ p(x,y)\in R \mapsto p(x+cyt,y+c^2t^2) + (t^3) \in R_3=R[[t]]/(t^3).$$
Since the coefficient of $t^3$ in $f(x+cyt,y+c^2t^2)$ is $0$, we deduce that $\delta$ is $I$-logarithmically $3$-integrable and so $\Ider_k(\log I;3)=\Ider_k(\log I;2)$.
A generic $I$-logarithmic $2$-integral of $\delta$ is determined by the $k$-algebra map
$$ p(x,y)\in R \mapsto p(x+cyt+dt^2,y-c^2t^2) + (t^3) \in R_3,$$
with $d\in R$, and a generic $I$-logarithmic $3$-integral of $\delta$ is determined by the $k$-algebra map
$$ p(x,y)\in R \mapsto p(x+cyt+dt^2+e t^3,y+c^2t^2) + (t^4) \in R_4,$$
with $d,e\in R$. The coefficient of $t^4$ in
$f(x+cyt+dt^2+e t^3,y+c^2t^2)$
is $d^2 + yc^4$, and so, the following conditions are equivalent:
\begin{enumerate}
\item[(a)] $\delta$ is $I$-logarithmically $4$-integrable.
\item[(b)] There is a $d\in R$ such that $d^2 + yc^4\in J$.
\end{enumerate}

The proof of the following lemma is easy:

\begin{lema} \label{lema:ej-nolineal} The set
$\Gamma:=\{c\in R\ |\ \exists d\in R,\ d^2 + yc^4\in J\}$
is the ideal generated by $x,y$.
\end{lema}

As a consequence of the lemma we deduce that $\{x\dpar{x},y^2\dpar{x},f\dpar{y}\}$ is a system of generators of
$\Ider_k(\log I;4)$. But $y^2\dpar{x}$ is $I$-logarithmically ($\infty$-)integrable after Proposition \ref{prop:Jder-in-Ider}, and so
\begin{eqnarray*}
&\Der_k(A)= \langle \overline{\dpar{x}}\rangle \varsupsetneq \Ider_k(A;2) = \langle \overline{x\dpar{x}}, \overline{y\dpar{x}}\rangle = \Ider_k(A;3) \varsupsetneq &\\
&\Ider_k(A;4) = \langle \overline{x\dpar{x}}, \overline{y^2\dpar{x}}\rangle = \Ider_k(A;5) = \cdots =\Ider_k(A;\infty).
\end{eqnarray*}
In particular, we have
\begin{eqnarray*}
&  \ann_A \left( \Der_k(A)/\Ider_k(A;2)\right) = ( \overline{x},\overline{y}) =\sqrt{\overline{J}} \varsupsetneq &\\ &  \ann_A \left( \Der_k(A)/\Ider_k(A;\infty)\right) = (\overline{x},\overline{y}^2) \varsupsetneq \overline{J}=(\overline{x}^2,\overline{y}^2).
\end{eqnarray*}
\medskip

Let us now compute the case $p=3$. The Jacobian ideal of $f$ is $J=(x,f)=(x,y^3)$. In a
similar way to the preceding case, we obtain that:
\begin{enumerate}
\item[-)]  $\Der_k(\log I) = \langle f\dpar{x},\dpar{y}\rangle$.
\item[-)] Since $2$ is invertible in $k$ we have $\Der_k(\log I) =\Ider_k(\log I;2)$.
\item[-)] $\Ider_k(\log I;3)=\langle x\dpar{y},y\dpar{y},f\dpar{x}\rangle$.
\item[-)] $\Ider_k(\log I;3) = \Ider_k(\log I;\infty)$.
\item[-)] $\Der_k(A)= \langle \overline{\dpar{y}}\rangle =\Ider_k(A;2)  \varsupsetneq \Ider_k(A;3) = \langle \overline{x\dpar{y}}, \overline{y\dpar{y}}\rangle = \Ider_k(A;4) = \cdots =\Ider_k(A;\infty)$ and $\ann_A \left( \Der_k(A)/\Ider_k(A;\infty)\right) = (\overline{x},\overline{y}) =\sqrt{J_{A/k}}$.
\end{enumerate}

Let us notice that for the cusp in characteristics $\neq 2,3$ we can apply
Proposition \ref{prop:isol-sing-general} and obtain that any derivation is
integrable.

\subsection{The cusp $x^2+y^3$ over the integers} \label{subsect:cusp-Z}

Assume that $k=\ZZ$ and  $f=x^2+y^3 \in R=\ZZ[x,y]$. Let $I=(f)$ and $A=\ZZ[x,y]/I$.
The Jacobian ideal of $f$ is $J=(2x,3y^2,f)=(2x,3y^2,x^2,y^3)$.
The $I$-logarithmic derivations of $R$ are generated by
$\delta_1=3x\dpar{x}+2y\dpar{y}$, $\delta_2=3y^2\dpar{x}-2x\dpar{y}$, $f\dpar{x}$ and $f\dpar{y}$.
The first derivation $\delta_1$ is the Euler vector field for the weights $w(x)=3, w(y)=2$.
As in \ref{subsect:cusp-2-3}, $\delta_1$ is $I$-logarithmically integrable. For the second derivation $\delta_2$, we apply Proposition \ref{prop:Jder-in-Ider} and we deduce that it is also $I$-logarithmically integrable. So this is an example of a non-smooth $\ZZ$-algebra $A$ for which any derivation is integrable.

\subsection{The cusp $3x^2+2y^3$ over the integers} \label{subsect:cusp2-Z}

Assume that $k=\ZZ$ and $f=3x^2+2y^3 \in R=\ZZ[x,y]$. Let $I=(f)$ and $A=\ZZ[x,y]/I$.
The Jacobian ideal of $f$ is $J=(6x,6y^2,f)=(6x,6y^2,3x^2,2y^3)$.
The $I$-logarithmic derivations of $R$ are generated by
$\delta_1=3x\dpar{x}+2y\dpar{y}$ and $\delta_2=-y^2\dpar{x}+x\dpar{y}$, which in
fact form a basis (we can say that ``$f$ is a free divisor'' of $R$). As in
\ref{subsect:cusp-2-3}, $\delta_1$ is the Euler vector field for the weights
$w(x)=3, w(y)=2$ and so it is $I$-logarithmically integrable.
\medskip

Let us study the integrability of $a\delta_2$, $a\in R$. The coefficient of $t^2$ in
$f(x-ay^2t,y+axt)$ is $a^2(3y^4+6x^2 y)$. Since $6x^2\in J$, this coefficient
belongs to $J$ if and only if $3a^2y^4 \in J$, i.e. $a^2\in J:3y^4$.

\begin{lema} \label{lema:ej-nolineal-2}
\begin{enumerate}
\item[(a)] $J:3y^4 = (2,x^2)$.
\item[(b)] $\{a\in R\ |\ a^2 \in (2,x^2)\} = (2,x)$.
\end{enumerate}
\end{lema}

\begin{corolario} The $R$-module $\Ider_\ZZ(\log I;2)$ is generated by $\{\delta_1,2\delta_2,x\delta_2\}$ and so
$ \ann_A \left( \Der_\ZZ(A)/\Ider_\ZZ(A;2)\right) = (2,x)$.
\end{corolario}

Let us study the $3$-integrability of
$$(2b+cx)\delta_2= -y^2(2b+cx)\dpar{x}+(2b+cx)x\dpar{y},\quad b,c\in R.$$
Let us write $a=2b+cx$. The coefficient of $t^2$ in $f(x-y^2(2b+cx)t,y+(2b+cx)xt)$
is $A(2y^3)+B(3x^2)$ with $A=6b(b+cx)y, B=c^2y^4 + 2a^2y$, which can be expressed as
$$ (A-B)xf'_x+(A-B)yf'_y+(3B-2A)f.$$
Hence, the coefficient of $t^2$ in
$$f(x-y^2(2b+cx)t+(B-A)xt^2,y+(2b+cx)xt+(B-A)yt^2)$$
is $(3B-2A)f$ and the reduction $\mod t^3$ of the $\ZZ$-algebra map
$$\Psi^{(2)}: p(x,y)\in R \mapsto p(x-y^2(2b+cx)t+(B-A)xt^2,y+(2b+cx)xt+(B-A)yt^2) \in R[[t]]$$
is $I$-logarithmic and gives rise to a $I$-logarithmic $2$-integral of $a\delta_2$.
So, the reduction $\mod t^3$ of the $\ZZ$-algebra map $\Psi_g^{(2)}:R\to R[[t]]$
given by
\begin{eqnarray*}
x & \mapsto &x-y^2(2b+cx)t+[(B-A)x+3dx-ey^2]t^2,\\
y & \mapsto &y+(2b+cx)xt+[(B-A)y+2dy+ex]t^2
\end{eqnarray*}
 is the associated map to a generic $I$-logarithmic $2$-integral of $a\delta_2$.
Moreover, the coefficient of $t^2$ in $\Psi_g^{(2)}(f)$ is $(3B-2A+6d)f$.
\medskip

The coefficient of $t^3$ in $\Psi_g^{(2)}(f)$ is
 $6 x^{2} y^{6} c^{3}+12 {x} y^{6} {b} c^{2}+12 x^{4} y^{3} c^{3}+36 x^{3} y^{3} {b} c^{2}+2 x^{6} c^{3}+36 x^{2} y^{3} b^{2} {c}+12 x^{5} {b} c^{2}+24 {x} y^{3} b^{3}+24 x^{4} b^{2} {c}+6 {x} y^{4} {c} {e}+16 x^{3} b^{3}+6 x^{2} y^{2} {c} {d}+12 y^{4} {b} {e}+12 x^{3} {y} {c} {e}+12 {x} y^{2} {b} {d}+24 x^{2} {y} {b} {e}$,
and it belongs to $J$ if and only if
$$ 2x^6c^3+16x^3b^3\in J \Leftrightarrow x^3c^3+8b^3 \in (J:2x^3). $$

\begin{lema} \label{lema:ej-nolineal-3} With the above notations, the following assertions hold:
\begin{enumerate}
\item[(a)] $J:2x^3 = (3,y^3)$.
\item[(b)] $ x^3c^3+8b^3 \in (J:2x^3)\Leftrightarrow a^3\in (J:2x^3) \Leftrightarrow
a\in (3,y)$.
\end{enumerate}
\end{lema}

\begin{corolario} The $I$-logarithmic derivation  $a\delta_2$ is $I$-logarithmically
$3$-integrable  if and only if $a\in (2,x)\cap (3,y) = (6,3x,2y,xy)$, and so the
$R$-module $\Ider_\ZZ(\log I;3)$ is generated by
$\{\delta_1,6\delta_2,3x\delta_2,2y\delta_2,xy\delta_2\}$ and
\begin{eqnarray*}
&\ann_A \left( \Der_\ZZ(A)/\Ider_\ZZ(A;3)\right) = (2,\overline{x})\cap (3,\overline{y}),&\\
& \ann_A \left( \Der_\ZZ(A;2)/\Ider_\ZZ(A;3)\right) =(3,\overline{y}) .
\end{eqnarray*}
\end{corolario}

The following lemma cannot be deduced directly from Proposition \ref{prop:Jder-in-Ider}. Its proof proceeds by induction and it is left up to the reader.

\begin{lema} \label{lema:tecnico-Z-cusp} Let $a\in (2,x)\cap (3,y)$. There are sequences $a_i, b_i\in R$, $i\geq 2$,
such that the $\ZZ$-algebra map
$$\Psi: p(x,y)\in R \mapsto p\left(x-ay^2t+\sum_{i=2}^\infty a_i t^i,y+axt+\sum_{i=2}^\infty b_i t^i\right) \in R[[t]]$$
is $I$-logarithmic, i.e. $\Psi(f) \in R[[t]]f$.
\end{lema}

\begin{corolario} We have
$$ \Ider_\ZZ(A;3) = \Ider_\ZZ(A;4)= \cdots = \Ider_\ZZ(A),$$
and so
$$ \ann_A \left( \Der_\ZZ(A)/\Ider_\ZZ(A)\right) = (2,\overline{x})\cap (3,\overline{y}) \supsetneq \sqrt{J_{A/\ZZ}} = (3\overline{x},2\overline{y}).$$
\end{corolario}

The following two examples have been proposed by Herwig Hauser.

\subsection{The surface $x_3^2+x_1 (x_1+x_2)^2=0$ in characteristic $2$}

Let $k$ be a field of characteristic $2$,  $f=x_3^2+x_1 (x_1+x_2)^2 \in
R=k[x_1,x_2,x_3]$, $I=(f)$ and $A=R/I$. The Jacobian ideal is $J=(\ell^2,f)=
(\ell^2,x_3^2)$ with $\ell=x_1+x_2$, 
and $\sqrt{J}=(\ell,x_3)$. A system of generators of $\Der_k(\log I)$ mod.
$f\Der_k(R)$ is $\{\dpar{2}, \dpar{3}\}$.
\medskip

\begin{lema} Let $\alpha,\beta\in R$ and $\delta=\alpha \dpar{2} + \beta \dpar{3}$.
The following conditions are equivalent:
\begin{enumerate}
\item[(a)] $\delta$ is $I$-logarithmically $2$-integrable.
\item[(b)] $x_1\alpha^2+\beta^2\in J$.
\end{enumerate}
\end{lema}

\begin{lema} The module $\{(\alpha,\beta)\in R^2\ |\ x_1\alpha^2+\beta^2\in
J\}$ is generated by $(x_3,0), (\ell,0),(0,x_3), (0,\ell)$.
\end{lema}

\begin{corolario} A system of generators of $\Ider_k(\log I;2)$ $\mod f\Der_k(R)$ is
$\{x_3 \dpar{2}, \ell \dpar{2}, x_3 \dpar{3}, \ell \dpar{3}\}$.
\end{corolario}

\begin{proposicion} $\Ider_k(A;2)=\Ider_k(A)$.
\end{proposicion}

\begin{prueba} We need to prove that $x_3 \dpar{2}, \ell \dpar{2}, x_3 \dpar{3}, \ell
\dpar{3}$ are $I$-logarithmically integrable. \smallskip

The derivation $x_3\dpar{3}$ is the Euler vector field for the weights
$w(x_1)=w(x_2)=2, w(x_3)=3$. From \ref{nume:traves} we deduce that $x_3\dpar{3}$ is
$I$-logarithmically integrable. \smallskip

The derivation $\ell \dpar{3}$ is $I$-logarithmically integrable since
$f(x_1+t^2,x_2+t^2,x_3 + \ell t) =\cdots= f\in R[t]\subset R[[t]]$ and so a
$I$-logarithmic integral of $\ell \dpar{3}$ is given by the $k$-algebra map $R \to
R[[t]]$ determined by $$ x_1 \mapsto x_1+t^2,\quad x_2 \mapsto x_2+t^2,\quad
x_3\mapsto x_3 + \ell t.$$ \smallskip

For the derivation $x_3 \dpar{2}$ let us write $W(t) = \frac{x_1^2 t^2}{1-x_1 t^2}
\in (t^2)R[[t]]$ and consider the homomorphism of $k$-algebras $\Psi: R \to R[[t]]$
given by:
$$x_1 \mapsto x_1+ W(t),\quad
x_2  \mapsto x_2 + x_3 t + W(t),\quad x_3 \mapsto x_3.$$ We have $\Psi(f)=f(x_1+
W,x_2 + x_3 t + W,x_3)=\cdots= \left(\frac{1}{1-x_1 t^2}\right) f$ and so $\Psi$
gives rise to a $I$-logarithmic integral of $x_3 \dpar{2}$. \smallskip

For the derivation $\ell \dpar{2}$ let us write $V(t) = \frac{x_1 t^2}{1-t^2} \in
(t^2) R[[t]]$ and consider the homomorphism of $k$-algebras $\Psi: R \to R[[t]]$
given by:
$$x_1 \mapsto x_1+ V(t),\quad x_2  \mapsto x_2 + \ell t + V(t),\quad
x_3 \mapsto x_3.$$ We have $\Psi(f)=f(x_1+ V,x_2 + \ell t + V,x_3)=\cdots = f$ and
so $\Psi$ gives rise to a $I$-logarithmic integral of $\ell \dpar{2}$.
\end{prueba}

 In this example the descending chain of modules of integrable
derivations stabilizes from $N=2$:
$$ \Der_k(A)=\Ider_k(A;1) \supset \Ider_k(A;2) =\Ider_k(A;3)=\cdots=
\Ider_k(A;\infty)$$ and
$$\ann_A \left( \Der_k(A)/\Ider_k(A;\infty)\right) = (\ell,x_3) = \sqrt{J}/I.$$

\subsection{The surface $x_3^2+x_1 x_2 (x_1+x_2)^2=0$ in characteristic $2$}

Let $k$ be a field of characteristic $2$,  $f=x_3^2+x_1 x_2 (x_1+x_2)^2 \in
R=k[x_1,x_2,x_3]$, $I=(f)$ and $A=R/I$. The Jacobian ideal is
$J=(x_2\ell^2,x_1\ell^2,f)= (x_2\ell^2,x_1\ell^2,x_3^2)$ with $\ell=x_1+x_2$.
It is clear that $\sqrt{J}=(\ell,x_3)$. The module $\Der_k(\log I)$ is generated
mod. $f\Der_k(R)$ by $\dpar{3}$, $\varepsilon = x_1\dpar{1} + x_2\dpar{2}$ and $\eta
= x_1^2 \ell^2\dpar{1}+x_3^2 \dpar{2}$ ($\dpar{3}(f)=\varepsilon(f)=0,
\eta(f)=x_1\ell^2f$). Since $\varepsilon$ is the Euler vector field for the weights
$w(x_1)=w(x_2)=1, w(x_3)=2$, we deduce from \ref{nume:traves} that $\varepsilon$ is
$I$-logarithmically integrable. From Proposition \ref{prop:Jder-in-Ider} we also
deduce that $\eta$ is $I$-logarithmically integrable.
\medskip

To find a system of generators of $\Ider_k(\log I;2)$ we need the conditions on
$a\in R$ which guarantee that $a\dpar{3}$ is $I$-logarithmically $2$-integrable. The
coefficient of $t^2$ in $f(x_1,x_2,x_3+at)=f+a^2t^2$ is $a^2$, and so $a\dpar{3}$ is
$I$-logarithmically $2$-integrable if and only if $a^2 \in J$.

\begin{lema} $\{a\in R | a^2 \in J\} = (x_3,x_1\ell,x_2\ell)$.
\end{lema}

\begin{corolario} A system of generators of $\Ider_k(\log I;2)$ mod. $f\Der_k(R)$ is
$\{x_3\dpar{3},x_1\ell\dpar{3},x_2\ell\dpar{3},\varepsilon,\eta\}$. In particular we
have
$$ \ann_A \left( \Der_k(A)/\Ider_k(A;2)\right) = (\overline{x_3},\overline{x_2} \overline{\ell},\overline{x_1}\overline{\ell)}.$$
\end{corolario}

The following lemma is a very particular case of a general result.

\begin{lema} Any Hasse--Schmidt derivation $E\in\HS_k(A;2)$ is $3$-integrable.
\end{lema}

\begin{prueba} Since $3$ is invertible in $k$, we can
consider the differential operator $E_3= E_1 E_2 - \frac{1}{3}E_1^3$ and check that
$(\Id_A,E_1,E_2,E_3)$ is a Hasse--Schmidt derivation.
\end{prueba}

As a consequence of the above lemma we have $\Ider_k(A;2)=\Ider_k(A;3)$.
\smallskip

Let us see the conditions for $a\dpar{3}$, with $a=\alpha x_3 +\beta x_1\ell+\gamma
x_2\ell$, $\alpha,\beta,\gamma\in R$, to be $I$-logarithmically $4$-integrable. The
algebra map associated with a general $I$-logarithmic $3$-integral of $a\dpar{3}$ is
$\Psi^{(3)}:R \to R_3$ given by:
\begin{eqnarray*}
x_1 & \mapsto &x_1+(\alpha^2 x_1 +\gamma^2 x_2+B_1x_1+C_1x_1^2\ell^2) t^2+(B_2x_1+C_2x_1^2\ell^2)t^3,\\
x_2 & \mapsto &x_2 + (\beta^2 x_1+B_1 x_2+C_1 x_3^2)t^2+(B_2 x_2+C_2 x_3^2)t^3,\\
x_3 &\mapsto &x_3+(\alpha x_3 +\beta x_1\ell+\gamma x_2\ell)t + A_1 t^2+A_2 t^3
\end{eqnarray*}
with $A_2,B_2,C_2\in R$, and let $\Psi_0^{(4)}:R\to R_4$ be the obvious lifting of
$\Psi^{(3)}$. The coefficient $\mod J$ of $t^4$ in the expression of
$\Psi_0^{(4)}(f)$, is $x_1 x_2^3(\alpha+\beta+\gamma)^4+A_1^2$. So, we have proved
the following lemma.

\begin{lema}  With the above notations, the following assertions are equivalent:
\begin{enumerate}
\item[(a)] The logarithmic derivation $a\dpar{3}$,
with $a=\alpha x_3 +\beta x_1\ell+\gamma x_2\ell$, is $I$-logarithmically
$4$-integrable.
\item[(b)] There is $A_1\in R$ such that $x_1 x_2^3(\alpha+\beta+\gamma)^4+A_1^2\in
J$, or, equivalently, $x_1 x_2^3(\alpha+\beta+\gamma)^4 \in J + R^2$.
\end{enumerate}
\end{lema}

\begin{lema} \label{lema:ej-nolineal-4} We have
$\{\varphi\in R\ |\ x_1 x_2^3 \varphi^4 \in J + R^2\}=(x_3,\ell)$.
\end{lema}

\begin{prueba} Let us write $\mathfrak{A}=\{\varphi\in R\ |\ x_1 x_2^3 \varphi^4 \in J +
R^2\}$. It is clear that $x_3,\ell\in \mathfrak{A}$, since $x_3^4\in J$ and $x_1
x_2^3 \ell^4\in J$. Let $\varphi$ be an element in $\mathfrak{A}$ and let us write
$\varphi = q x_3 +\varphi_1(x_1,x_2)$, with $q\in R$ and $\varphi_1(x_1,x_2)\in
\mathfrak{A}$. We have
$$x_1 x_2^3 \varphi_1^4 = U(x_1,x_2) x_1 \ell^2 + V(x_1,x_2) x_2 \ell^2 +
P(x_1,x_2)^2.$$ By taking derivatives with respect to $x_1$ we obtain $x_2^3
\varphi_1^4 = U'_{x_1} x_1 \ell^2 +U \ell^2 + V'_{x_1} x_2\ell^2$ and so $\ell$
divides $\varphi_1$. We conclude that $\mathfrak{A} = (x_3,\ell)$.
\end{prueba}

As a consequence of the above lemma and the fact that $(x_3,\ell)$ is a prime ideal,
the condition $x_1 x_2^3(\alpha+\beta+\gamma)^4 \in J + R^2$ is equivalent to
$\alpha+\beta+\gamma \in (x_3,\ell)$, i.e. to $\alpha = \alpha_1 x_3 + \alpha_2 \ell
+ \beta +\gamma$ and so $ a = \cdots = \alpha_1 x_3^2 + \alpha_2 x_3\ell +
\beta(x_3+x_1\ell) +\gamma(x_3+x_2\ell)$. We conclude with the following corollary.

\begin{corolario} A system of generators of $\Ider_k(\log I;4)$ $\mod f\Der_k(R)$
is
$\{x_3^2\dpar{3},x_3\ell\dpar{3},(x_3+x_1\ell)\dpar{3},(x_3+x_2\ell)\dpar{3},\varepsilon,\eta\}$.
In particular we have
\begin{eqnarray*}
& \ann_A \left( \Der_k(A)/\Ider_k(A;2)\right) = (\overline{x_3},\overline{x_2} \overline{\ell},\overline{x_1}\overline{\ell}),&\\
& \ann_A \left( \Der_k(A)/\Ider_k(A;4)\right) =
(\overline{x_3}^2,\overline{x_3}\overline{\ell},\overline{x_3}+\overline{x_2}\overline{\ell},\overline{x_3}+\overline{x_1}\overline{\ell}),&\\
& \ann_A \left( \Ider_k(A;2)/\Ider_k(A;4)\right) = (\overline{x_3},\overline{\ell})&
\end{eqnarray*} and all the inclusions
$$ J_{A/k} \subset (\overline{x_3}^2,\overline{x_3}\overline{\ell},\overline{x_3}+\overline{x_2}\overline{\ell},\overline{x_3}+\overline{x_1}\overline{\ell})
\subset (\overline{x_3},\overline{x_2}
\overline{\ell},\overline{x_1}\overline{\ell}) \subset
(\overline{x_3},\overline{\ell})=\sqrt{J_{A/k}}$$ are strict.
\end{corolario}

From Proposition \ref{prop:Jder-in-Ider} we deduce that $x_3^2\dpar{3}$ is
$I$-logarithmically integrable.

\begin{lema} \label{lema:tecnico-xld3}
The derivation $x_3\ell\dpar{3}$ is $I$-logarithmically integrable.
\end{lema}

\begin{prueba} Let us write $\delta=x_3\ell\dpar{3}$ and $D=(x_3\ell)\sbullet
\Delta^{(3)}$. We have $\Phi_D(f)= f + (x_3\ell)^2 t^2$ and $(x_3\ell)^2 =
f_{x_1}'f_{x_2}'+\ell^2 f = x_1 x_2 \ell^4 + \ell^2 f$. Let us also write
$S=k[x_1,x_2]$ and $\mathfrak{b}=(f_{x_1}',f_{x_2}')=(x_2 \ell^2,x_1 \ell^2)\subset
S$.

We are going to construct inductively a sequence of differential operators
$E^m_m\in\mathfrak{b}\diff_{S/k}$, $m\geq 1$, with $E^1_1=0$, $E^2_2(f)= x_1 x_2
\ell^4$, $E^m_m(f)=0$ for all $m\geq 3$ and such that
$(\Id,E^1_1,E^2_2,E^3_3,\dots)$ is a Hasse-Schmidt derivation of length $\infty$.
\smallskip

For $m=2$, let us take $E^2_2 =f_{x_2}'\dpar{1}$.
\smallskip

Assume that we have already found a Hasse--Schmidt derivation
$E^m=(\Id,E^1_1,\dots,E^m_m)\in\HS_k(S;m)$ with the required properties. Let us
consider $F^m=\varepsilon(E^m)\in \HS_k(S;\infty)$. From Proposition
\ref{prop:formulon}, 2) we deduce that $F^m_{m+1} \in \mathfrak{b}^2 \diff_{S/k}$
and so $F^m_{m+1}(f)\in \mathfrak{b}^2$. Hence, there are $\alpha,\beta \in
\mathfrak{b}$ such that $F^m_{m+1}(f) = \alpha f_{x_1}'+\beta f_{x_2}'$ and
consequently we can take $E^{m+1}_{m+1}=F^m_{m+1} - (\alpha \dpar{1}+\beta\dpar{2})$
\smallskip

Once the Hasse--Schmidt derivation $E=(\Id,0,E^2_2,E^3_3,\dots)\in \HS_k(S;\infty)$
has been constructed, we extend it in the obvious way to the ring $R$ (we keep the
same name $E$ for the extension). We have $\Phi_{D\pcirc E}(f) =
\widetilde{\Phi}_D\left( \Phi_E(f)\right) = \widetilde{\Phi}_D\left( f + x_1 x_2
\ell^4 t^2 \right) = \Phi_{D}(f) + \Phi_{D}(x_1 x_2 \ell^4)t^2 = f + (x_3\ell)^2 t^2
+ x_1 x_2 \ell^4t^2 = (1+\ell^2 t^2)f$
and so $D\pcirc E$ is a $I$-logarithmic integral of $\delta$.
\end{prueba}

The proof of the following lemma is due to M. M\'erida.

\begin{lema} 
The derivations $(x_3+x_1\ell)\dpar{3}$ and $(x_3+x_2\ell)\dpar{3}$ are $I$-logarithmically integrable.
\end{lema}

\begin{prueba} By symmetry, it is enough to consider the case $(x_3+x_1\ell)\dpar{3}$, for which the logarithmic integrability is a consequence of the fact that the map $\Psi: R \to R[[t]]$ given by:
\begin{eqnarray*}
x_1 & \mapsto &x_1+x_1V,\\
x_2 & \mapsto &x_2+x_1V,\\
x_3 &\mapsto &x_3+(x_3+x_1\ell)t+x_3V,
\end{eqnarray*}
with $\displaystyle V = \sum_{i=1}^\infty t^{2^i}$, is $I$-logarithmic. Namely, since $t^2=V^2+V$, we have
\begin{eqnarray*}
& f(x_1+x_1V,x_2+x_1V,x_3+(x_3+x_1\ell)t+x_3V) =&\\
&(x_3+(x_3+x_1\ell)t+x_3V)^2 + (x_1+x_1V)(x_2+x_1V)\ell^2=&\\
&x_3^2 +(x_3^2+x_1^2 \ell^2) t^2 +x_3^2 V^2 + (x_1 x_2 + x_1^2 V + x_1 x_2 V + x_1^2 V^2)\ell^2=&\\
& x_3^2 +(x_3^2+x_1^2 \ell^2) t^2 +x_3^2 V^2 + (x_1 x_2 + x_1^2 t^2 + x_1 x_2 V)\ell^2 = &\\
& x_3^2 +x_3^2 t^2 +x_3^2 V^2 + (x_1 x_2 +  x_1 x_2 V)\ell^2 = f + x_3^2 t^2 +x_3^2 V^2 + x_1 x_2 V\ell^2= &\\
& f + x_3^2 V + x_1 x_2 V\ell^2 = (1+V) f.&
\end{eqnarray*}
\end{prueba}

\begin{corolario}  $\Ider_k(A;4)=\Ider_k(A)$.
\end{corolario}

\subsection{Some questions}

\begin{question} \label{cuestion:1} Assume that $R=k[x_1,\dots,x_d]$, $S\subset R$
is a multiplicative set and $A=S^{-1}R$ or $A=k[[x_1,\dots,x_d]]$. Let $I\subset A$
be an ideal, $m\geq 1$ an integer, $D\in\HS_k(\log I;m)$ and
$E=\overline{D}\in\HS_k(A/I;m)$. Let us consider the following properties:
\begin{enumerate}
\item[(a)] $D$ is $I$-logarithmically $n$-integrable for all integers $n\geq m$ (or equivalently, $E$ is $n$-integrable for all integers $n\geq m$).
\item[(b)] $D$ is $I$-logarithmically $\infty$-integrable (or equivalently $E$ is $\infty$-integrable).
\end{enumerate}
Under which hypotheses on $k$ and on $I$ are properties (a) and (b) equivalent for
any $D\in\HS_k(\log I;m)$? Are they equivalent if $k$ is a field or the ring of
integers and $I$ is arbitrary?
\medskip

Notice that this question is the same as asking whether the inclusion in
\ref{eq:intersection-iderlog}  (or in \ref{eq:intersection-ider} for $m=1$) is an equality
or not.
\end{question}

\begin{question} The proofs of propositions \ref{prop:surjec-IHS-log-localiz-poly} and
\ref{prop:surjec-IHS-localiz-fp} do not work for $m=\infty$ and, presumably, these
propositions are not true for $m=\infty$ without additional finiteness hypotheses on
$k$. Let us notice that if the maps in Proposition \ref{prop:surjec-IHS-localiz-fp}
are surjective for $m=\infty$, then the localization conjecture for the
Hasse--Schmidt algebra stated in \cite{traves-2003} is true.
\end{question}

\begin{question} For any finitely presented $k$-algebra $A$,
find an algorithm for deciding whether \underline{a given} $\delta\in\Der_k(A)$ is
$m$-integrable or not.
\end{question}

\begin{question}  For any finitely presented $k$-algebra $A$,
find an algorithm to obtain a system of generators of $\Ider_k(A;m)$, $m\geq 2$.
\end{question}

\begin{question} Assume that the base ring $k$ is a field of positive characteristic
or $\ZZ$, or perhaps a more general noetherian ring, and $A$ a finitely generated
$k$-algebra. Is there an integer $n\geq 1$ such that $\Ider_k(A;n)=
\Ider_k(A;\infty)$? Or at least, is the descending chain of $A$-modules $
\Ider_k(A;1) \supset  \Ider_k(A;2) \supset  \Ider_k(A;3) \supset \cdots$ stationary?
\end{question}

\begin{question} Assume that the base ring $k$ is a field of positive characteristic
or $\ZZ$, or perhaps a more general noetherian ring. Is there an integer $m\gg 1$,
possibly depending on $d$ and $e$ or other numerical invariants, such that
$$ \Ider_k(A;m)=\Der_k(A)\quad \Rightarrow\quad \Ider_k(A)=\Der_k(A)$$
for every quotient ring $A=k[x_1,\dots,x_d]/I$ with $\dim A= e$?
\end{question}

\begin{question} Assume that the base ring $k$ is a field of positive characteristic
or $\ZZ$, or perhaps a more general noetherian ring, $A$ a local noetherian
$k$-algebra and $\delta:A\to A$ a $k$-derivation. Under which hypotheses the
$m$-integrability of $\widehat{\delta}:\widehat{A} \to \widehat{A}$ implies the
$m$-integrability of $\delta$?
\end{question}

{\small \noindent \href{http://departamento.us.es/da/}{Departamento de \'{A}lgebra} \&\
Instituto de
Matem\'aticas (\href{http://www.imus.us.es}{IMUS})\\
Facultad de Matem\'aticas, Universidad de Sevilla\\
P.O. Box 1160, 41080
 Sevilla, Spain}. \\
{\small {\it E-mail}\ : narvaez@algebra.us.es
 }

\end{document}